\documentclass{amsart}
\usepackage{amssymb}
\usepackage{mathrsfs}
\usepackage{stmaryrd}
\usepackage{bbm}
\usepackage{oldgerm}
\usepackage[francais]{babel}
\usepackage[all]{xy}

\newtheorem{teo}[subsection]{Th\'eor\`eme}
\newtheorem{prop}[subsection]{Proposition}
\newtheorem{cor}[subsection]{Corollaire}
\newtheorem{lem}[subsection]{Lemme}
\newtheorem{conj}[subsection]{Conjecture}

\theoremstyle{definition}

\newtheorem{defi}[subsection]{D\'efinition}
\newtheorem{rema}[subsection]{Remarque}

\numberwithin{equation}{subsection}

\newcommand{\mQ}{{\mathbb Q}}
\newcommand{\mN}{{\mathbb N}}
\newcommand{\mZ}{{\mathbb Z}}
\newcommand{\mF}{{\mathbb F}}
\newcommand{\mA}{{\mathbb A}}

\newcommand{\mG}{{\mathbb G}}

\newcommand{\mK}{{\mathbb K}}
\newcommand{\bG}{{\bf G}}
\newcommand{\id}{{\rm id}}
\newcommand{\Ob}{{\rm ob\ }}
\newcommand{\pr}{{\rm pr}}
\newcommand{\Gr}{{\rm Gr}}
\newcommand{\gr}{{\rm gr}}
\newcommand{\iso}{{\rm iso}}
\newcommand{\Aut}{{\rm Aut}}

\newcommand{\fil}{{\rm fil}}
\newcommand{\ord}{{\rm ord}}
\newcommand{\res}{{\rm res}}

\newcommand{\m}{{\rm m}}

\newcommand{\fC}{{\mathfrak C}}

\newcommand{\fE}{{\mathfrak E}}
\newcommand{\fF}{{\mathfrak F}}
\newcommand{\fX}{{\mathfrak X}}
\newcommand{\fY}{{\mathfrak Y}}

\newcommand{\fm}{{\mathfrak m}}
\newcommand{\rig}{{\rm rig}}

\newcommand{\Spec}{{\rm Spec}}
\newcommand{\Proj}{{\rm Proj}}
\newcommand{\et}{{\rm \Acute{e}t}}

\newcommand{\sw}{{\rm sw}}
\newcommand{\rsw}{{\rm rsw}}

\newcommand{\rD}{{\rm D}}

\newcommand{\bD}{{\bf D}}

\newcommand{\bP}{{\bf P}}

\newcommand{\bT}{{\bf T}}
\newcommand{\bV}{{\bf V}}
\newcommand{\Hom}{{\rm Hom}}
\newcommand{\cHom}{{\mathscr Hom}}

\newcommand{\cF}{{\mathscr F}}
\newcommand{\cG}{{\mathscr G}}
\newcommand{\cD}{{\mathscr D}}

\newcommand{\cL}{{\mathscr L}}

\newcommand{\cJ}{{\mathscr J}}
\newcommand{\cI}{{\mathscr I}}
\newcommand{\co}{{\mathscr O}}
\newcommand{\cS}{{\mathscr S}}
\newcommand{\cK}{{\mathscr K}}
\newcommand{\cH}{{\mathscr H}}

\newcommand{\cN}{{\mathscr N}}

\newcommand{\hA}{{\widehat A}}

\newcommand{\rB}{{\rm B}}
\newcommand{\rC}{{\rm C}}
\newcommand{\rF}{{\rm F}}
\newcommand{\rH}{{\rm H}}
\newcommand{\rN}{{\rm N}}
\newcommand{\rR}{{\rm R}}
\newcommand{\rV}{{\rm V}}
\newcommand{\rW}{{\rm W}}
\newcommand{\rZ}{{\rm Z}}

\newcommand{\oa}{\overline{a}}

\newcommand{\ok}{\overline{k}}
\newcommand{\oq}{\overline{q}}
\newcommand{\os}{\overline{s}}

\newcommand{\oK}{\overline{K}}
\newcommand{\oR}{\overline{R}}
\newcommand{\oF}{\overline{F}}
\newcommand{\oG}{\overline{G}}

\newcommand{\oS}{\overline{S}}

\newcommand{\obT}{\overline{\bT}}
\newcommand{\uk}{{\underline{k}}}

\newcommand{\tx}{\widetilde{x}}
\newcommand{\ty}{\widetilde{y}}
\newcommand{\tI}{\widetilde{I}}
\newcommand{\tJ}{\widetilde{J}}

\newcommand{\ta}{\widetilde{a}}

\newcommand{\ttC}{{\tt C}}
\newcommand{\ttS}{{\tt S}}

\newcommand{\ttb}{{\tt b}}
\newcommand{\ttc}{{\tt c}}

\newcommand{\tth}{{\tt h}}

\newcommand{\ttt}{{\tt t}}
\newcommand{\ttu}{{\tt u}}
\newcommand{\ttv}{{\tt v}}
\newcommand{\ttx}{{\tt x}}
\newcommand{\tty}{{\tt y}}
\newcommand{\ttty}{\widetilde{\tt y}}
\newcommand{\ctf}{{\tt ctf}}
\newcommand{\vc}{{\tt SS}}

\newcommand{\oeta}{{\overline{\eta}}}

\newcommand{\oTheta}{{\overline{\Theta}}}

\newcommand{\trC}{\widetilde{\ttC}}

\begin{document}

\title[Analyse micro-locale $\ell$-adique en caract\'eristique $p>0$]
{Analyse micro-locale $\ell$-adique en caract\'eristique $p>0$\\
Le cas d'un trait}
\author{Ahmed Abbes, Takeshi Saito}
\address{A.A. CNRS UMR 7539, LAGA, Institut Galil\'ee, Universit\'e Paris 13,
93430 Villetaneuse, France}
\email{abbes@math.univ-paris13.fr}
\address{T.S. Department of Mathematical Sciences, 
University of Tokyo, Tokyo 153-8914, Japan}
\email{t-saito@ms.u-tokyo.ac.jp}

\begin{abstract}
Nous d\'eveloppons, pour un faisceau \'etale $\ell$-adique
sur un trait complet de caract\'eris\-tique $p>0$, la notion de vari\'et\'e
caract\'eristique. Notre approche, inspir\'ee de l'analyse micro-locale 
de Kashiwara et Schapira, est un pendant faisceautique de notre th\'eorie
de ramification  des corps locaux \`a corps r\'esiduel quelconque.
Nous pr\'esentons la principale propri\'et\'e que devrait satisfaire
la  vari\'et\'e caract\'eristique (conjecture de l'isog\'enie), et 
nous la d\'emontrons pour les faisceaux de rang un sans 
restriction sur le trait, ou inconditionnellement sur le faisceau 
si le corps r\'esiduel du trait est parfait.\\[2mm]
{\sc Abstract.}\ We develop, for an $\ell$-adic \'etale sheaf on a 
complete trait of characteristic $p>0$, the notion of characteristic
variety. Our approach, inspired by the microlocal analysis of Kashiwara
and Schapira, is a complement to our ramification theory for local
fields with general residue fields. We formulate the main property 
that should be satisfied by the characteristic variety (the isogeny 
conjecture), and prove it for rank one sheaves unconditionally 
on the trait, or unconditionally on the sheaf if the residue field of 
the trait is perfect.  
\end{abstract}

\keywords{Corps locaux, ramification, faisceaux \'etales $\ell$-adiques,
vari\'et\'e caract\'eristique, caract\`eres d'Artin-Schreier-Witt, 
conducteurs de Swan de Kato.}

\subjclass[2000]{primary 11S15,
secondary 14F20}

\maketitle

\section{Introduction}

Ce travail s'inscrit dans un projet consacr\'e \`a g\'en\'eraliser, pour 
les faisceaux \'etales $\ell$-adiques sur une vari\'et\'e de caract\'eritique
$p>0$ (avec $\ell\not= p$), la notion emprunt\'ee \`a la th\'eorie des 
$\cD$-modules de {\em vari\'et\'e caract\'eristique}. Pour un $\cD$-module
coh\'erent sur une vari\'et\'e analytique complexe, 
la th\'eorie {\em d'analyse micro-locale} 
de Kashiwara et Schapira permet de reconstruire 
la vari\'et\'e caract\'eristique \`a partir du complexe de ses solutions 
holomorphes (\cite{ksc} 11.3.3). C'est donc tout naturellement que nous nous 
en inspirons pour notre projet. 
Il y a eu une premi\`ere tentative dans cette direction due \`a
Verdier \cite{verdier}; elle n'a pas abouti car le foncteur de 
sp\'ecialisation de Verdier tue la ramification sauvage. 
 
Dans cet article, nous d\'eveloppons la variante locale en codimension 
un du projet.  Consid\'erons un corps parfait $k$
de caract\'eristique $p>0$ et un anneau de valuation  
discr\`ete complet $R$ qui est une $k$-alg\`ebre, dont le corps r\'esiduel 
est de type fini sur $k$. 
On pose $S=\Spec(R)$, $\eta$ (resp. $s$) son point 
g\'en\'erique (resp. ferm\'e), $\oeta$ un point g\'eom\'etrique g\'en\'erique
de $S$ et $\os$ le point g\'eom\'etrique localis\'e en $s$ correspondant. 
Soit $D$ un diviseur effectif de $S$. Suivant Raynaud, on introduit
dans \eqref{dilat16} un ouvert $(S\times_kS)_{(D)}$ de l'\'eclatement 
de $S\times_kS$ le long de $D$ diagonalement plong\'e, 
d\'efini par une condition de valuation. 
On l'appelle {\em dilatation de $D$} et on le 
consid\`ere comme un $S$-sch\'ema par la seconde projection. 
Sa fibre g\'en\'erique est isomorphe \`a $S\times_k\eta$
et sa fibre au dessus de $D$ est isomorphe au fibr\'e vectoriel 
$\bT_D=\bV(\Omega^1_R\otimes_R\co_D(D))$;
on pose $\bT^\ttt_D$ le fibr\'e vectoriel dual sur $D$. 
On peut voir dans cette construction une analogie 
avec la d\'eformation au c\^one normal de Grothendieck, qui est aussi 
un exemple de dilatation \eqref{def-cn}. 

Soit $\Lambda$ un anneau commutatif noeth\'erien  
annul\'e par un entier inversible dans $k$. 
On d\'efinit, pour les faisceaux \'etales en $\Lambda$-modules 
sur $S\times_k\eta$,
le foncteur $\nu_D$ de {\em sp\'ecialisation le long de $D$}  
comme \'etant le foncteur 
des cycles proches de $(S\times_kS)_{(D)}$, et le foncteur 
$\mu_D$ de {\em micro-localisation le long de $D$} comme \'etant
la transform\'ee de Fourier-Deligne de $\nu_D$ relativement \`a un 
caract\`ere additif non-trivial.

L'anneau $R$ appara\^{\i}t comme l'anneau local compl\'et\'e d'un $k$-sch\'ema
lisse $X$ en un point de codimension $1$. 
Il est beaucoup plus ais\'e de travailler avec le $S$-sch\'ema 
$(X\times_kS)_{(D)}$ obtenu en dilatant $X\times_kS$, plut\^ot
que $(S\times_kS)_{(D)}$. Pour d\'emontrer nos r\'esultats, nous nous y
ramenons par alg\'ebrisation. 
Pour les applications g\'eom\'etriques que nous avons en vue, 
un tel sch\'ema $X$ s'impose naturellement. Donc on aurait pu 
se limiter \`a ce cadre, mais cela aurait l'inconv\'enient 
de rendre la th\'eorie locale d\'ependante d'un choix superflu. 
Il est aussi possible de dilater le sch\'ema formel compl\'et\'e 
de $S\times_kS$ le long de $D$ diagonalement plong\'e, 
ce qui reviendrait \`a compl\'eter 
$(S\times_kS)_{(D)}$ le long de sa fibre sp\'eciale. 
Nous avons choisi de ne pas poursuivre cet approche ici car elle requi\`ere 
le d\'eveloppement du formalisme \'etale pour le sch\'ema 
formel ainsi obtenu et pour sa fibre rigide au sens de Raynaud.
Nous reviendrons \`a cet aspect dans un travail ult\'erieur.

Soient $\cF$ un faisceau \'etale
constructible en $\Lambda$-modules plats sur $\eta$, 
$\cF_!$ son extension par $0$ \`a $S$. 
Notre id\'ee principale, inspir\'ee de la th\'eorie de Kashiwara et Schapira,
est de consid\'erer le support $\ttC_D(\cF)$ dans 
$\obT_D^\ttt=\bT_D^\ttt\times_D\os$
du complexe $\mu_D(\cF_!\boxtimes\Lambda)$, 
o\`u $D$ est un {\em divisieur critique} pour $\cF$. 
Pour d\'efinir ces diviseurs, on utilise la th\'eorie de ramification 
pour les corps locaux \`a corps r\'esiduel quelconque \cite{as},
dont ce travail est un pendant faisceautique. Nous avons d\'efini dans 
{\em loc. cit.} une filtration d\'ecroissante, exhaustive et s\'epar\'ee
$(G^a)_{a\in \mQ_{\geq 0}}$ de $G={\rm Gal}(\oeta/\eta)$ 
par des sous-groupes ferm\'es distingu\'es. 
Les {\em pentes critiques} de $\cF$ sont les nombres rationnels $r\geq 0$
tel que $(\cF_\oeta)^{G^r}\subsetneq (\cF_\oeta)^{G^{r+}}$. 
Ces derniers n'\'etant pas entiers en g\'en\'eral, nous sommes amen\'es 
\`a modifier notre construction en rempla\c{c}ant le second facteur 
de $S\times_kS$ par une extension finie. 

Supposons que $\cF$ ait une unique pente critique $r>1$, 
et pour simplifier l'introduction que $r$ soit entier, 
ce qui d\'etermine un diviseur effectif $D$ de $S$. 
Nous conjecturons que {\em $\ttC_D(\cF)$ est un ensemble 
fini de points du fibr\'e vectoriel $\obT^\ttt_D$
ne contenant pas l'origine} \eqref{ciso1}. 
Concr\`etement, cela permet d'associer \`a $\cF$ un nombre fini 
de formes diff\'erentielles ``tordues'' non-nulles. 
En homog\'en\'eisant $\ttC_D(\cF)$, on obtient la 
{\em vari\'et\'e caract\'eristique} de $\cF$ \eqref{def-vc}. 
Cette d\'efinition ne depend pas du caract\`ere additif choisi.

Notre r\'esultat principal 
est la d\'emonstration de la conjecture ci-dessus en rang un. 
Pour ce faire, on compare notre approche 
\`a la th\'eorie de ramification de Kato 
pour les caract\`eres \cite{kato1}. 
On montre \eqref{ntp1} que la pente critique d'un faisceau de rang 
un correspond 
\`a son conducteur de Swan, et le support $\ttC_D(\cF)$ \`a son conducteur 
de Swan raffin\'e. On retrouve en fait la variante 
de la th\'eorie de Kato modifi\'ee par Matsuda \cite{mats}; 
mais on d\'eveloppe aussi une variante logarithmique 
qui correspond exactement \`a la th\'eorie de Kato en rang un \eqref{tp1}. 
Enfin, on d\'emontre la conjecture inconditionnellement 
si le corps r\'esiduel de $R$ est $k$, donc parfait \eqref{parf2}.

L'article se compose de deux parties qui ont chacune ses propres notations. 
La premi\`ere partie d\'eveloppe la notion g\'eom\'etrique 
de dilatation. La seconde partie est consacr\'ee \`a l'analyse micro-locale; 
les constructions principales sont introduites dans les sections 
\ref{spml} et \ref{lspml}; 
les principaux r\'esultats sont \'enonc\'es dans la section \ref{enonces}
et d\'emontr\'es dans les sections \ref{Dtp1} et \ref{Dtp2}. 
La th\'eorie de ramification 
des caract\`eres d'Artin-Schreier-Witt de Kato est reprise enti\`erement 
dans les sections \ref{thkato} et \ref{thkato-mats}.
Enfin, on calcule dans la derni\`ere section les invariants de ramification 
d\'efinis dans cet article en termes d'invariants plus classiques
sous l'hypoth\`ese que le corps r\'esiduel de $R$ est $k$.

\part{Dilatations}

\section{Rel\`evement d'un morphisme aux sch\'emas \'eclat\'es} 

\subsection{}\label{ann-grad}
Soient $\ttS$, $\ttS'$ deux anneaux gradu\'es \`a degr\'es positifs, 
$\theta\colon \ttS'\rightarrow \ttS$ un homomorphisme d'anneaux gradu\'es. 
Suivant EGA II 2.8.1, on d\'esigne par $\bG(\theta)$
la partie ouverte de $X=\Proj(\ttS)$ r\'eunion des $\rD_+(\theta(f'))$, o\`u
$f'$ parcourt l'ensemble des \'el\'ements homog\`enes de $\ttS'_+$. 
En vertu de {\em loc. cit.} 2.8.2 et 2.8.3, 
il existe un morphisme affine 
canonique du sch\'ema induit $\bG(\theta)$ dans $\Proj(\ttS')$, dit 
associ\'e \`a $\theta$. 
Soient $\ttS''$  un troisi\`eme anneau gradu\'e \`a degr\'es positifs, 
$\theta'\colon \ttS''\rightarrow \ttS'$ un homomorphisme d'anneaux gradu\'es, 
et posons $\theta''=\theta\circ \theta'$. 
D'apr\`es {\em loc. cit.}  2.8.4, on a $\bG(\theta'')\subset \bG(\theta)$,
et si $\vartheta$, $\vartheta'$ et $\vartheta''$ sont les morphismes 
associ\'es \`a $\theta$, $\theta'$ et $\theta''$, 
on a $\vartheta(\bG(\theta''))\subset \bG(\theta')$ et 
$\vartheta''=\vartheta'\circ (\vartheta|\bG(\theta''))$.

\subsection{}\label{alg-grad}
Soient $Y$ un sch\'ema, $\cS$, $\cS'$ deux $\co_Y$-alg\`ebres gradu\'ees
quasi-coh\'erentes \`a degr\'es positifs; posons $X=\Proj(\cS)$, 
$X'=\Proj(\cS')$, et soient $p$, $p'$ les morphismes structuraux de $X$ et 
$X'$ dans $Y$. Soit $\theta\colon \cS'\rightarrow \cS$ un $\co_Y$-homomorphisme
d'alg\`ebres gradu\'ees.  
Pour tout ouvert affine $U$ de $Y$, on pose $\ttS_U=\Gamma(U,\cS)$, 
$\ttS'_U=\Gamma(U,\cS)$; l'homomorphisme $\theta$ d\'efinit 
un homomorphisme $\theta_U\colon \ttS'_U\rightarrow \ttS_U$ 
de $A_U$-alg\`ebres gradu\'ees;
en posant $A_U=\Gamma(U,\co_Y)$. Il lui correspond dans $p^{-1}(U)$ un
ensemble ouvert $\bG(\theta_U)$ et un morphisme 
$\vartheta_U\colon \bG(\theta_U)
\rightarrow p'^{-1}(U)$. D'apr\`es EGA II 3.5.1, il existe une partie
ouverte $\bG(\theta)$ de $X$ et un $Y$-morphisme affine 
$\vartheta\colon \bG(\theta)\rightarrow X'$ dit associ\'e \`a $\theta$,
tels que, pour tout ouvert affine $U$ de $Y$, 
$\bG(\theta)\cap p^{-1}(U)=\bG(\theta_U)$ et $\vartheta_U$ soit
la restriction de $\vartheta$ au dessus de $U$. 

\subsection{}\label{dilat1}
Soient $X,Y$ deux sch\'emas, $f\colon Y\rightarrow X$ un morphisme, 
$D$ (resp. $E$) un sous-sch\'ema ferm\'e de $X$ (resp. $Y$) 
d\'efini par un id\'eal quasi-coh\'erent 
$\cI$ (resp. $\cJ$) de $\co_X$ (resp. $\co_Y$), 
$i\colon D\rightarrow X$ et $j\colon E\rightarrow Y$ les injections 
canoniques. On d\'esigne par $X'$ l'\'eclatement de $D$
dans $X$ et par $Y'$ l'\'eclatement de $E$ dans $Y$. 
On suppose que  $f\circ j$ soit major\'e par $i$; il 
revient au m\^eme de dire que $f^*(\cI)\co_Y\subset \cJ$. 
Par suite, on a un homomorphisme de $\co_Y$-alg\`ebres gradu\'ees 
\begin{equation}
\theta\colon f^*(\oplus_{n\geq 0}\cI^n)\rightarrow \oplus_{n\geq 0}\cJ^n
\end{equation}
qui d\'etermine par \eqref{alg-grad} un ensemble ouvert 
$Y'_{(D)}=\bG(\theta)$ de $Y'$ 
et un morphisme $\vartheta\colon Y'_{(D)}\rightarrow Y\times_XX'$.

\begin{defi}[Raynaud \cite{blr} 3.2]\label{dilat2}
On appelle $Y'_{(D)}$ la {\em dilatation} de $E$ par rapport \`a $D$,  
$\vartheta\colon Y'_{(D)}\rightarrow Y\times_XX'$ le
{\em morphisme canonique}, et le compos\'e $f'\colon Y'_{(D)}\rightarrow X'$ 
de $\vartheta$ et de la projection canonique le {\em rel\`evement canonique} 
de $f$. 
\end{defi}  

\subsection{}\label{dilat25}
Sous les hypoth\`eses de \ref{dilat1}, si $E$ est l'image r\'eciproque 
de $D$, ce qui revient \`a dire que $f^*(\cI)\co_Y=\cJ$, 
alors $Y'_{(D)}=Y'$ et le morphisme canonique $\vartheta$ 
est une immersion ferm\'e (EGA II 3.6.2). 
Si de plus $E$ est un diviseur de Cartier de $Y$, alors $Y'$ 
est isomorphe \`a $Y$ 
et le rel\`evement canonique $f'\colon Y\rightarrow X'$ traduit 
la propri\'et\'e universelle des \'eclatements.

\subsection{}\label{dilat3}
Conservons les hypoth\`eses de \ref{dilat1}, de plus soient 
$Z$ un troisi\`eme sch\'ema, $g\colon Z\rightarrow Y$
un morphisme, $F$ un sous-sch\'ema ferm\'e de $Z$,
$k\colon F\rightarrow Z$ l'injection canonique,  
$Z'$ l'\'eclatement de $F$ dans $Z$. 
On suppose que $g\circ k$ soit major\'e par $j$. On peut alors consid\'erer 
les dilatations de $F$ par rapport \`a $E$ ou \`a $D$, 
et les morphismes canoniques 
$\vartheta'\colon Z'_{(E)}\rightarrow Z\times_YY'$ et 
$\vartheta''\colon Z'_{(D)}\rightarrow Z\times_XX'$. Il r\'esulte 
de \ref{ann-grad} qu'on a $Z'_{(D)}\subset Z'_{(E)}$, $\vartheta'(Z'_{(D)})
\subset Z\times_YY'_{(D)}$ et $\vartheta''=(Z\times_Y\vartheta)\circ 
(\vartheta'|Z'_{(D)})$.

\begin{prop}\label{dilat4}
Les hypoth\`eses \'etant celles de \eqref{dilat1}, de plus soit 
$g\colon Z\rightarrow Y$ un morphisme tel que l'id\'eal $\cJ\co_Z$ 
soit inversible et $\cI\co_Z=\cJ\co_Z$. 
Alors  le rel\`evement canonique $g'\colon Z\rightarrow Y'$ de $g$
se factorise \`a travers l'ouvert $Y'_{(D)}$ de $Y'$. 
\end{prop}  

En effet, $Z$ s'identifie aux dilatations de 
$g^{-1}(E)=g^{-1}(f^{-1}(D))$ par rapport \`a $E$ ou \`a $D$ \eqref{dilat25}. 
Donc la proposition r\'esulte de la fonctorialit\'e des dilatations
\eqref{dilat3}. 

\begin{prop}\label{dilat5}
Sous les hypoth\`eses de \eqref{dilat1}, 
on a l'\'egalit\'e des id\'eaux 
\begin{equation}\label{dilat6}
(\cI\co_{Y'})|Y'_{(D)}=(\cJ\co_{Y'})|Y'_{(D)}
\end{equation}
et $Y'_{(D)}$ est l'ouvert maximal de $Y'$ o\`u cette relation est satisfaite. 
\end{prop}

Comme $\cI\co_{Y'}\subset \cJ\co_{Y'}$ et $\cJ\co_{Y'}$ est un id\'eal 
inversible, $(\cI\co_{Y'})\otimes (\cJ\co_{Y'})^{-1}$ est un id\'eal 
quasi-coh\'erent de $\co_{Y'}$; il d\'efinit un sous-sch\'ema
ferm\'e de $Y'$ dont le compl\'ementaire $U$ est l'ouvert
maximal de $Y'$ tel que $(\cI\co_{Y'})|U=(\cJ\co_{Y'})|U$. 
Il r\'esulte de \ref{dilat4} que $U$ est contenu dans $Y'_{(D)}$. 
Pour \'etablir l'\'egalit\'e $U=Y'_{(D)}$, il suffit de montrer la 
relation \eqref{dilat6}. La question \'etant locale sur $X$ et $Y$, 
on peut supposer que $X=\Spec(A)$ et  $Y=\Spec(B)$ soient affines;
donc $f$ est d\'efini par un homomorphisme $u\colon A\rightarrow B$. 
Soient $I$ (resp. $J$) l'id\'eal de $A$ (resp. $B$) tel que $\cI=\tI$
(resp. $\cJ=\tJ$). On pose $\ttS=\oplus_{n\geq 0}J^n$, de sorte que 
$Y'=\Proj(\ttS)$. Pour $a\in J$, si on note $a'$ l'\'el\'ement $a$ vu comme 
un \'el\'ement homog\`ene de degr\'e un de $\ttS$,  
l'image inverse de $a$ engendre l'id\'eal $J\co_{Y'}$ sur l'ouvert $\rD_+(a')$.
Donc l'inclusion $I\co_{Y'}\subset J\co_{Y'}$ induit une \'egalit\'e 
sur les ouverts $\rD_+(u(a)')$ pour $a\in I$.  
Mais on a $Y'_{(D)}=\cup_{a\in I}\rD_+(u(a)')$, 
ce qui entra\^{\i}ne \eqref{dilat6} et ach\`eve la preuve.

\begin{cor}\label{dilat7}
Sous les hypoth\`eses de \eqref{dilat3}, on a 
$Z'_{(D)}=g^{\prime -1}(Y'_{(D)})$ o\`u $g'\colon Z'_{(E)}\rightarrow Y'$
est le rel\`evement canonique de $g$. 
\end{cor}

On pose $U=g^{\prime -1}(Y'_{(D)})$ et on note
$\cK$ l'id\'eal de $\co_{Z}$ qui d\'efinit $F$. On sait \eqref{dilat3} que 
$Z'_{(D)}\subset U\subset Z'_{(E)}$. D'autre part, on a 
$(\cI\co_{Y'})|Y'_{(D)}=(\cJ\co_{Y'})|Y'_{(D)}$ 
et $(\cJ\co_{Z'})|Z'_{(E)}=(\cK\co_{Z'})|Z'_{(E)}$.  
Donc $(\cI\co_{Z'})|U=(\cK\co_{Z'})|U$, ce qui 
entra\^{\i}ne que $U\subset Z'_{(D)}$ \eqref{dilat4}.

\begin{cor}\label{dilat8}
Les hypoth\`eses \'etant celles de \eqref{dilat1}, de plus notons
$E'$ le diviseur exceptionnel sur $Y'$, c'est \`a dire l'image r\'eciroque
de $E$ sur $Y'$, $E'_{(D)}$ sa restriction \`a $Y'_{(D)}$. 
Alors le diagramme 
\[
\xymatrix{
{E'_{(D)}}\ar[r]\ar[d]&{Y'_{(D)}}\ar[d]\\
{D}\ar[r]&{X}}
\]
est cart\'esien.
\end{cor}
Cela r\'esulte de la relation \eqref{dilat6}.

\subsection{}\label{def-cn}
Soient $k$ un corps, $Y$ un $k$-sch\'ema, 
$X$ un sous-sch\'ema ferm\'e de $Y$, $i\colon X\rightarrow Y$ l'injection 
canonique. Consid\'erons le diagramme commutatif 
\[
\xymatrix{X\ar[r]^j\ar[d]&{\mA^1_Y}\ar[d]^\pi\\
{\Spec(k)}\ar[r]_{s_0}&{\mA^1_k}}
\]
o\`u $\pi$ est la projection canonique, $s_0$ est la section nulle de $\mA^1_k$
et $j$ est le compos\'e de $i$ et de la section nulle de $\mA^1_Y$.  
La dilatation de $X$ dans $\mA^1_Y$ par rapport \`a l'origine de 
$\mA^1_k$ n'est autre que l'espace total de la d\'eformation au c\^one 
normal; sa restriction au dessus de $\mG_{\m,k}$ est isomorphe \`a $\mG_{\m,Y}$
et sa fibre au dessus de l'origine de $\mA^1_k$ est isomorphe au c\^one 
normal de $X$ dans $Y$.

\section{Dilatation le long d'une section}

\subsection{}\label{dilat10} 
Soient $X,Y$ deux sch\'emas, $f\colon Y\rightarrow X$ un morphisme s\'epar\'e, 
$g\colon X\rightarrow Y$ une section de $f$; 
donc $g$ est une immersion ferm\'ee. 
Soient $D$ un sous-sch\'ema ferm\'e de $X$,
$i\colon D\rightarrow X$ l'injection canonique.
L'immersion compos\'ee $g\circ i\colon D\rightarrow Y$
d\'etermine un sous-sch\'ema ferm\'e de $Y$ qu'on notera simplement $g(D)$.   
On d\'esigne par $X_{[D]}$ l'\'eclatement de $D$ dans $X$, par $Y_{[D]}$ 
l'\'eclatement de $g(D)$ dans $Y$ et par $Y_{(D)}$ la dilatation de 
$g(D)$ par rapport \`a $D$ \eqref{dilat2}. On dira que $Y_{(D)}$ est la
{\em dilatation de $Y$ le long de la section $g$ d'\'epaisseur $D$}. 
On notera $E_{[D]}$ (resp. $E_{(D)}$) le diviseur exceptionnel sur 
$Y_{[D]}$ (resp. $Y_{(D)}$).
D'apr\`es \ref{dilat5} et \ref{dilat8}, on a un diagramme 
commutatif canonique \`a carr\'es cart\'esiens 
\begin{equation}\label{dilat105} 
\xymatrix{
{E_{(D)}}\ar[r]\ar[d]&{Y_{(D)}}\ar[d]&{Y\times_XU}\ar[l]\ar[d]\\
D\ar[r]&X&U\ar[l]}
\end{equation}

\subsection{}\label{dilat11} 
Les hypoth\`eses \'etant celles de \eqref{dilat10}, 
de plus soient $u\colon X'\rightarrow X$ un morphisme,
$D'=u^{-1}(D)$ l'image r\'eciproque de $D$ sur $X'$,   
$Y'=X'\times_XY$, $f'\colon Y'\rightarrow X'$
et $v\colon Y'\rightarrow Y$ les projections canoniques,
$g'\colon X'\rightarrow Y'$ la section d\'eduite de $g$
par le changement de base $u$. On a alors un diagramme commutatif 
\begin{equation}\label{dilat12}
\xymatrix{
{Y'_{[D']}}\ar@{}[rd]|{\Box}\ar[d]_{v_{[D]}}&{Y'_{(D')}}\ar[l]
\ar[d]^{v_{(D)}}\ar[r]
&{X'_{[D']}}\ar[d]^{u_{[D]}}\ar[r]&X'\ar[d]^{u}\\
{Y_{[D]}}&{Y_{(D)}}\ar[l]\ar[r]&{X_{[D]}}\ar[r]&X}
\end{equation}
o\`u $v_{[D]}, v_{(D)}$ et $u_{[D]}$ sont les morphismes canoniques.  
En effet, $g'(D')$ est l'image inverse de $g(D)$ par $v$, ce qui justifie 
la d\'efinition de $v_{[D]}$. Compte tenu de \ref{dilat5}, $Y'_{(D')}$ 
s'identifie aussi \`a la dilatation de $g'(D')$ par rapport \`a $D$; 
la fl\`eche $v_{(D)}$ et la commutativit\'e
de \eqref{dilat12} r\'esultent alors de la fonctorialit\'e des dilatations
\eqref{dilat3}. 
Le carr\'e de gauche est cart\'esien par \ref{dilat7}. 
Si de plus $u$ est plat, les trois carr\'es sont cart\'esiens;
ceci est \'evident pour le carr\'e de droite; donc la fl\`eche canonique 
$Y_{(D)}\times_{X}X'\rightarrow  Y'$ satisfait aux conditions 
de la proposition \ref{dilat4}, et induit un morphisme 
$Y_{(D)}\times_{X}X'\rightarrow  Y'_{(D')}$;
on laissera au lecteur le soin d'en d\'eduire que le carr\'e central
est cart\'esien. 

\subsection{}\label{fp}
Conservons les hypoth\`eses de \eqref{dilat10}, de plus soient 
$Z$ un sch\'ema, $p\colon Z\rightarrow Y$ un morphisme s\'epar\'e, 
$h\colon X\rightarrow Z$ un morphisme tel que $p\circ h=g$; donc $h$
est une section de $f\circ p$. Il r\'esulte de \ref{dilat3}
qu'on a un diagramme commutatif 
\begin{equation}\label{fp1}
\xymatrix{
{E^Z_{(D)}}\ar@{}[rd]|{\Box}\ar[r]\ar[d]_{p_E}&
{Z_{(D)}}\ar[r]\ar[d]^{p_{(D)}}&Z\ar[d]^p\\
{E^Y_{(D)}}\ar[r]&{Y_{(D)}}\ar[r]&Y}
\end{equation}
o\`u $E^Y_{(D)}$ (resp. $E^Z_{(D)}$) est le diviseur exceptionnel sur 
$Y_{(D)}$ (resp. $Z_{(D)}$); de plus le carr\'e de gauche est cart\'esien 
par \ref{dilat8}. Si le diagramme 
\begin{equation}\label{fp2}
\xymatrix{
D\ar[r]^{h\circ i}\ar@{=}[d]&Z\ar[d]^p\\
D\ar[r]_{g\circ i}&Y}
\end{equation}
est cart\'esien, \eqref{fp1} s'ins\`ere dans un diagramme 
commutatif \`a carr\'e central cart\'esien (cf. \ref{dilat7}) 
\begin{equation}\label{fp3}
\xymatrix{
{E^Z_{(D)}}\ar@{}[rd]|{\Box}\ar[r]\ar[d]_{p_E}&
{Z_{(D)}}\ar@{}|\Box[rd]\ar[r]\ar[d]|{p_{(D)}}&
{Z_{[D]}}\ar[r]\ar[d]^{p_{[D]}}&Z\ar[d]^p\\
{E^Y_{(D)}}\ar[r]&{Y_{(D)}}\ar[r]&{Y_{[D]}}\ar[r]&Y}
\end{equation}
Si de plus $p$ est plat, les trois carr\'es ci-dessus sont cart\'esiens
et $p_E$ est un isomorphisme.

\subsection{}\label{dilat125}
Conservons les hypoth\`eses de \eqref{dilat10}, supposons de plus 
que $D$ soit un diviseur de Cartier sur $X$ et posons 
$\co_D(D)=\co_X(D)\otimes_{\co_X}\co_D$.
On d\'esigne par $\cJ$ l'id\'eal de $\co_Y$ qui d\'efinit $g(D)$.
Le diagramme commutatif 
\[
\xymatrix{
D\ar[r]^i\ar@{=}[d]&X\ar[r]^{g}&Y\ar[d]^f\\
D\ar[rr]_i&&X}
\]
induit une suite exacte canonique de faisceaux conormaux
\begin{equation}\label{dilat135}
i^*(\cN_g)\rightarrow\cJ/\cJ^2\rightarrow\cN_i\rightarrow 0
\end{equation}
et un scindage $\cN_i\rightarrow \cJ/\cJ^2$. On en d\'eduit
un homomorphisme surjectif
\[
(\cN_g\otimes_{\co_X}\co_D(D))\oplus \co_D\rightarrow 
(\cJ/\cJ^2)\otimes_{\co_D}\co_D(D),
\]
et donc une immersion ferm\'ee $j$ de $E_{[D]}$ 
dans le fibr\'e projectif 
$\bP((\cN_g\otimes_{\co_X}\co_D(D))\oplus \co_D)$ au dessus de $D$
\footnote{Nous adoptons les conventions de EGA II 
pour les fibr\'es projectifs (4.1.1) et les fibr\'es vectoriels (1.7.8).}. 
Ce dernier est canoniquement isomorphe \`a la fermeture projective du
fibr\'e vectoriel $\bV(\cN_g\otimes_{\co_X}\co_D(D))$ avec un lieu 
\`a l'infini isomorphe \`a $\bP(\cN_g\otimes_{\co_X}\co_D(D))$.
Il r\'esulte de la d\'efinition \eqref{dilat1} que l'image inverse
de $\bV(\cN_g\otimes_{\co_X}\co_D(D))$ par $j$ est l'ouvert 
$E_{(D)}$. Identifiant $\cN_g$ avec $g^*(\Omega^1_{Y/X})$,
on obtient une immersion ferm\'ee au dessus de $D$
\begin{equation}\label{dilat14}
E_{(D)}\rightarrow \bV(g^*(\Omega^1_{Y/X})\otimes_{\co_X}\co_D(D)).
\end{equation}

\begin{prop}\label{dilat15}
Sous les hypoth\`eses de \eqref{dilat125}, 
si $g$ est une immersion r\'eguli\`ere, 
\eqref{dilat14} est un isomorphisme.
\end{prop}
En effet, $g\circ i$ est une immersion r\'eguli\`ere, 
la suite \eqref{dilat135} est exacte \`a gauche et 
l'immersion $j$ est un isomorphisme.

\subsection{}\label{dilat16}
Soient $k$ un corps, $X$ un $k$-sch\'ema, $Y$ un $k$-sch\'ema s\'epar\'e,
$f\colon X\rightarrow Y$ un $k$-morphisme, 
$D$ un divieur de Cartier effectif sur $X$, 
$U$ l'ouvert compl\'ementaire de $D$ dans $X$. 
On peut appliquer la construction \eqref{dilat10} relativement 
\`a la projection canonique $\pr_2\colon Y\times_kX\rightarrow X$ 
et au graphe $\gamma\colon X\rightarrow Y\times_kX$ de $f$;
on obtient un diagramme commutatif \`a carr\'es cart\'esiens
\begin{equation}\label{dilat164} 
\xymatrix{
{E_{(D)}}\ar[r]\ar[d]&{(Y\times_kX)_{(D)}}\ar[d]&
{Y\times_kU}\ar[l]\ar[d]&\\
D\ar[r]&X&U\ar[l]}
\end{equation}
et une $D$-immersion ferm\'ee 
$E_{(D)}\rightarrow \bV(f^*(\Omega^1_{Y/k})\otimes_{\co_X}\co_D(D))$
\eqref{dilat14}, 
qui est bijective si $Y$ est lisse sur $k$ \eqref{dilat15}.

\begin{lem}\label{dilat18}
Les hypoth\`eses \'etant celles de \eqref{dilat16}, supposons de plus 
qu'il existe un sous-sch\'ema ferm\'e $C$ de $Y$, 
que $F=f^{-1}(C)$ soit un diviseur de Cartier sur $X$ et que $D\geq 2 F$. 
Alors $C\times_kU$ est ferm\'e dans $(Y\times_kX)_{(D)}$. 
\end{lem} 
Le diagramme commutatif \`a carr\'es cart\'esiens canonique 
\[
\xymatrix{
F\ar@{=}[r]\ar[d]&F\ar[r]\ar[d]&C\ar[d]\\
D\ar[r]&X\ar[r]&Y}
\]
induit un diagramme cart\'esien 
\[
\xymatrix{
F\ar[r]\ar[d]&D\ar[d]\\
{C\times_kX}\ar[r]^{h}&{Y\times_kX}}
\]
Soit $h'\colon (C\times_kX)_{[F]}\rightarrow (Y\times_kX)_{[D]}$ 
la transform\'ee stricte de $h$; c'est une immersion ferm\'ee.
Il suffit de montrer que $(C\times_kX)_{[F]}\cap (Y\times_kX)_{(D)} 
=C\times_kU$, ou que le diviseur exceptionnel de $(C\times_kX)_{[F]}$
ne rencontre par $E_{(D)}$, autrement dit qu'il se plonge dans $E_{[D]}$
\`a travers le lieu \`a l'infini 
$\bP(f^*(\Omega^1_{Y/k})\otimes_{\co_X}\co_X(D))$ (cf. \ref{dilat125}).
On se r\'eduit \`a montrer que le morphisme canonique de faisceaux 
conormaux 
\[
\co_D(-D)\otimes_{\co_D}\co_F\rightarrow \cN_{F}(C\times_kX) 
\]
d\'eduit de la factorisation $F\rightarrow C\times_kD\rightarrow C\times_kX$
est nul. D'une part, cette derni\`ere s'ins\`ere dans la 
factorisation $F\rightarrow C\times_kF\rightarrow 
C\times_kD\rightarrow C\times_kX$. D'autre part, 
le morphisme canonique de faisceaux conormaux 
$\co_D(-D)\otimes_{\co_D}\co_F\rightarrow \co_F(-F)$ est nul \`a cause 
de l'in\'egalit\'e $D\geq 2F$.

\begin{lem}\label{dilat17}
Les hypoth\`eses \'etant celles de \eqref{dilat16}, soient de plus 
$T$ un $k$-sch\'ema, $f_1\colon X\rightarrow T$ et $f_2\colon T\rightarrow Y$
deux $k$-morphismes tels que $f=f_2\circ f_1$.   
Supposons que $f_2$ soit plat et que le diagramme 
\[
\xymatrix{
D\ar@{=}[d]\ar[r]^{f_1\circ i}&T\ar[d]^{f_2}\\
D\ar[r]_{f\circ i}&Y}
\]
soit cart\'esien. Alors on a un diagramme commutatif \`a carr\'es cart\'esiens
\[
\xymatrix{
{E^T_{(D)}}\ar[r]\ar[d]_{e_2}&{(T\times_kX)_{(D)}}\ar[r]\ar[d]&
{(T\times_kX)_{[D]}}\ar[d]\ar[r]&T\ar[d]^{f_2}\\
{E^Y_{(D)}}\ar[r]&{(Y\times_kX)_{(D)}}\ar[r]&
{(Y\times_kX)_{[D]}}\ar[r]&Y}
\]
o\`u on a not\'e $E^T_{(D)}$ et $E^Y_{(D)}$ les diviseurs exceptionnels;
de plus $e_2$ est un isomorphisme. 
\end{lem}

Cela r\'esulte de \ref{fp}.

\section{Produit fibr\'e logarithmique}

\subsection{}\label{plog1}
Soient $k$ un corps, 
$X$ (resp. $Y$) un $k$-sch\'ema, $X_0$ (resp. $Y_0$)
un sous-sch\'ema ferm\'e de $X$ (resp. $Y$). 
Consid\'erons le diagramme commutatif canonique
\[
\xymatrix{
{Y_0}\ar[d]&{Y_0\times_kX_0}\ar[l]\ar[r]\ar[d]&{X_0}\ar[d]\\
Y&{Y\times_kX}\ar[r]\ar[l]&X}
\]
On note $(Y\times_kX)'$ l'\'eclatement de $Y_0\times_kX_0$ dans $Y\times_kX$,
$(Y\times_kX)'_{(X_0)}$ (resp. $(Y\times_kX)'_{(Y_0)}$) la dilatation 
de $Y_0\times_kX_0$ dans $Y\times_kX$ par rapport \`a $X_0$
(resp. $Y_0$) et 
\begin{equation}\label{plog11}
(Y\times_kX)'_{(Y_0+X_0)}=(Y\times_kX)'_{(X_0)}\cap (Y\times_kX)'_{(Y_0)}.
\end{equation}

\subsection{}\label{plog15}
Conservons les hypoth\`eses de \eqref{plog1}, de plus soient 
$p\colon S\rightarrow X$, $q\colon T\rightarrow Y$ deux morphismes 
de sch\'emas, $S_0=p^{-1}(X_0)$ (resp. $T_0=q^{-1}(Y_0)$)  
l'image r\'eciproque de $X_0$ (resp. $Y_0$). On note $(T\times_kS)'$
l'\'eclatement de $T_0\times_kS_0$ dans $T\times_kS$ et 
$(T\times_kS)'_{(T_0+S_0)}$ l'ouvert d\'efini comme dans \eqref{plog11}.
Alors on a un diagramme commutatif 
\begin{equation}
\xymatrix{
{(T\times_kS)'_{(T_0+S_0)}}\ar[r]\ar[d]\ar@{}[rd]|{\Box}&
{(T\times_kS)'}\ar[r]\ar[d]&
{T\times_kS}\ar[d]^{q\times p}\\
{(Y\times_kX)'_{(Y_0+X_0)}}\ar[r]&{(Y\times_kX)'}\ar[r]&
{Y\times_kX}}
\end{equation}
o\`u le carr\'e de gauche est cart\'esien \eqref{dilat7}. 
Si $q$ et $p$ sont plats, les deux carr\'es sont cart\'esiens.

\subsection{}\label{plog3}
Conservons les hypoth\`eses de \eqref{plog1}, 
supposons de plus que $X_0$ et $Y_0$ soient des diviseurs de Cartier sur 
$X$ et $Y$ respectivement. Lorsqu'il n'y a aucun risque de confusion, 
on utilisera la notation $Y\times_k^{\log}X$ 
pour d\'esigner $(Y\times_kX)'_{(Y_0+X_0)}$, les diviseurs 
$X_0$ et $Y_0$ \'etant sous-entendus. 
Les immersions ferm\'ees $Y_0\times_kX\rightarrow Y\times_kX$ 
et $Y\times_kX_0\rightarrow Y\times_kX$ se rel\`event en des 
immersions ferm\'ees $Y_0\times_kX\rightarrow (Y\times_kX)'$ 
et $Y\times_kX_0\rightarrow (Y\times_kX)'$. On voit aussit\^ot 
que $Y\times_k^{\log}X$
est l'ouvert compl\'ementaire de $Y_0\times_kX$ et $Y\times_kX_0$ dans 
$(Y\times_kX)'$. 
On note $U$ (resp. $V$) l'ouvert 
compl\'ementaire de $X_0$ dans $X$ (resp. $Y_0$ dans $Y$). 
On a un diagramme commutatif \`a carr\'es cart\'esiens 
\[
\xymatrix{
V\ar[d]&{V\times_kU}\ar[r]\ar[l]\ar[d]&U\ar[d]\\
Y&{Y\times_k^{\log}X}\ar[l]^-(0.5){\pr_1}\ar[r]_-(0.5){\pr_2}&X}
\]
o\`u les fl\`eches horizontales sont les projections canoniques. 
Soient $f\colon X\rightarrow
Y$ un $k$-morphisme tel que $f^{-1}(Y_0)=X_0$,
$\gamma\colon X\rightarrow Y\times_kX$ son graphe. 
Comme $\gamma^{-1}(Y_0\times_kX_0)=X_0$, il r\'esulte de la propri\'et\'e 
universelle des dilatations \eqref{dilat4} 
que $\gamma$ se rel\`eve uniquement 
en une immersion $\gamma^{\log}\colon X\rightarrow Y\times_k^{\log}X$, 
dite {\em graphe logarithmique} de $f$. 

\begin{rema}\label{plog2}
{\rm (i)} Si dans \eqref{plog3}, $X=\Spec(A)$ et $Y=\Spec(B)$ sont affines
et $X_0$ (resp. $Y_0$) est d\'efini
par une \'equation $s\in A$ (resp. $t\in B$), alors 
$Y\times_k^{\log}X$ est le sch\'ema affine associ\'e \`a l'alg\`ebre 
\[
B\otimes_kA[\ttu,\ttu^{-1}]/(t\otimes 1-1\otimes s \cdot \ttu).
\]

{\rm (ii)} Si dans \eqref{plog3}, 
$X_0$ et $Y_0$ sont lisses sur $k$, on retrouve un cas particulier 
du produit logarithmique d\'efini par 
Kato et Saito (\cite{ks1} 4.2.4, \cite{ks2} 1.1.1).
\end{rema}

\begin{lem}\label{plog35}
Les hypoth\`eses \'etant celles de \eqref{plog3}, supposons de plus 
que $Y$ et $Y_0$ soient lisses sur $k$.  
Alors le morphisme $\pr_2\colon Y\times_k^{\log}X\rightarrow X$ 
est lisse et on a un isomorphisme canonique
\[
\Omega^1_{(Y\times_k^{\log}X)/X}\simeq \pr_1^*(\Omega^1_{Y/k}(\log Y_0)),
\]
o\`u $\Omega^1_{Y/k}(\log Y_0)$ est le faisceau des diff\'erentielles 
de $Y$ sur $k$ \`a p\^oles logarithmiques le long de $Y_0$.
\end{lem} 
Les questions \'etant locales, on peut se borner au cas o\`u 
$X=\Spec(A)$ et $Y=\Spec(B)$ soient affines 
et $X_0$ (resp. $Y_0$) soit d\'efini par une \'equation $s\in A$ (resp. 
$t\in B$). Alors $Y\times_k^{\log}X$ est le sch\'ema 
affine associ\'e \`a l'alg\`ebre 
$C=B\otimes_kA[\ttu,\ttu^{-1}]/(t\otimes 1-1\otimes s \cdot \ttu)$
\eqref{plog2}. Comme $B$ et $B/tB$ sont des $k$-alg\`ebres lisses, 
la suite 
\[
0\longrightarrow C\stackrel{\rho}{\longrightarrow} 
(\Omega^1_{B/k}\otimes_BC)\oplus Cd\ttu 
\longrightarrow \Omega^1_{C/A}\longrightarrow 0,
\]
o\`u $\rho(1)=dt\otimes 1-(1\otimes s)d\ttu$, 
est exacte et localement scind\'ee.
Le lemme s'ensuit par le crit\`ere jacobien.

\subsection{}\label{plog4}
Conservons les hypoth\`eses de \eqref{plog3}, supposons de plus 
que $Y$ soit s\'epar\'e sur $k$. 
Soient $f\colon X\rightarrow Y$ un $k$-morphisme 
tel que $f^{-1}(Y_0)=X_0$, $D$ un diviseur de Cartier effectif sur 
$X$ de m\^eme support que $X_0$. 
On peut appliquer la construction \eqref{dilat10}
relativement \`a la projection canonique 
$\pr_{2}\colon Y\times_k^{\log}X\rightarrow X$ (qui est s\'epar\'ee)
et au graphe logarithmique 
$\gamma^{\log}\colon X\rightarrow Y\times_k^{\log}X$ de $f$; 
on obtient un diagramme commutatif \`a carr\'es cart\'esiens
\begin{equation}\label{plog6} 
\xymatrix{
{E_{(D)}}\ar[r]\ar[d]&{(Y\times_k^{\log}X)_{(D)}}
\ar[d]&{V\times_kU}\ar[l]\ar[d]&\\
D\ar[r]&X&U\ar[l]}
\end{equation}
En vertu de \ref{dilat15} et \ref{plog35}, 
si $Y$ et $Y_0$ sont lisses sur $k$, on a un isomorphisme canonique
\begin{equation}\label{plog65}
E_{(D)}\stackrel{\sim}{\rightarrow} 
\bV(f^*(\Omega^1_{Y/k}(\log Y_0))\otimes_{\co_X}\co_D(D)).
\end{equation}

\begin{lem}\label{plog7}
Les hypoth\`eses \'etant celles de \eqref{plog4}, soient de plus 
$T$ un $k$-sch\'ema, $f_1\colon X\rightarrow T$ et $f_2\colon T\rightarrow Y$
deux $k$-morphismes tels que $f=f_2\circ f_1$.   
Supposons que $f_2$ soit plat et que le diagramme 
\[
\xymatrix{
D\ar@{=}[d]\ar[r]^{f_1\circ i}&T\ar[d]^{f_2}\\
D\ar[r]_{f\circ i}&Y}
\]
soit cart\'esien. Alors $T_0=f_2^{-1}(Y_0)$
est un diviseur de Cartier sur $T$, et 
on a un diagramme commutatif \`a carr\'es cart\'esiens
\[
\xymatrix{
{E^T_{(D)}}\ar[r]\ar[d]_{e_2}&{(T\times_k^{\log}X)_{(D)}}\ar[r]\ar[d]&
{(T\times_k^{\log}X)_{[D]}}\ar[d]\ar[r]&T\ar[d]^{f_2}\\
{E^Y_{(D)}}\ar[r]&{(Y\times_k^{\log}X)_{(D)}}\ar[r]&
{(Y\times^{\log}_kX)_{[D]}}\ar[r]&Y}
\]
o\`u on a not\'e $E^T_{(D)}$ et $E^Y_{(D)}$ les diviseurs exceptionnels;
de plus $e_2$ est un isomorphisme. 
\end{lem}

On note $p\colon T\times_k^{\log}X\rightarrow Y\times_k^{\log}X$
la fl\`eche d\'eduite de $f_2$ \eqref{plog15} et 
$\gamma_1^{\log}\colon X\rightarrow T\times_k^{\log}X$ (resp. 
$\gamma^{\log}\colon X\rightarrow Y\times_k^{\log}X$) le graphe logarithmique
de $f_1$ (resp. $f$). On a $\gamma^{\log}=p\circ \gamma_1^{\log}$ 
\eqref{dilat4}. Il r\'esulte de \ref{plog15} que le diagramme 
\[
\xymatrix{
D\ar[r]^-(0.5){\gamma_1\circ i}\ar@{=}[d]&
{T\times_k^{\log}X}\ar[r]\ar[d]^p&T\ar[d]^{f_2}\\
D\ar[r]_-(0.5){\gamma\circ i}&{Y\times_k^{\log}X}\ar[r]&Y}
\]
est commutatif et ses carr\'es sont cart\'esiens.
Le lemme r\'esulte alors de \ref{fp}.

\part{Analyse micro-locale $\ell$-adique}

\section{Notations et conventions}

\subsection{}
On d\'esigne par $k$ un corps parfait 
de caract\'eristique $p>0$. Pour toute $k$-alg\`ebre $A$, 
on note $\Omega^1_A=\Omega^1_{A/A^p}
=\Omega^1_{A/k}$ le module des $1$-diff\'erentielles absolues de $A$. 

\subsection{}\label{not1} 
On d\'esigne par $R$ un anneau de valuation 
discr\`ete complet qui est une $k$-alg\`ebre. 
On note $K$ le corps des fractions de $R$, $F$ son corps r\'esiduel 
qu'on suppose {\em de type fini sur $k$}, $\fm$ son id\'eal maximal, 
$\pi$ une uniformisante et $\ttv\colon K^\times \rightarrow \mZ$ la 
valuation de $K$ normalis\'ee par $\ttv(\pi)=1$. 
On pose $S=\Spec(R)$, $\eta=\Spec(K)$, $s=\Spec(F)$, 
$\oK$ une cl\^oture s\'eparable de $K$, 
$G$ le groupe de Galois de $\oK$ sur $K$, $\oR$ la cl\^oture 
int\'egrale de $R$ dans $\oK$, $\oF$ son corps r\'esiduel,
$\oS=\Spec(\oR)$, $\oeta=\Spec(\oK)$ et $\os=\Spec(\oF)$.

\subsection{}\label{diff1}
On notera que $R$ est non-canoniquement $k$-isomorphe \`a $F[[\pi]]$
(EGA IV 0.19.6.6). 
Comme $F$ est de type fini sur $k$, l'endomorphisme de Frobenius de $R$
est fini. Par cons\'equent, $\Omega^1_R$ est un $R$-module de type fini;
en particulier, il est complet et s\'epar\'e pour la topologie $\fm$-adique. 
D'autre part, $R$ est une $k$-alg\`ebre formellement lisse
pour la topologie $\fm$-adique (EGA IV 0.22.2.2).
Il r\'esulte alors que $\Omega^1_R$ est un $R$-module libre de type fini 
(EGA IV 0.20.4.10); on note $d$ son rang.

\subsection{}\label{diff3}
On pose $\Omega^1_R(\log)$ le $R$-module des $1$-diff\'erentielles 
logarithmiques absolues d\'efini par 
\begin{equation}
\Omega^1_R(\log)=(\Omega^1_R\oplus (R\otimes_{\mZ}K^\times))/\cF
\end{equation}
o\`u $\cF$ est le sous-$R$-module de 
$\Omega^1_R\oplus (R\otimes_{\mZ}K^\times)$
engendr\'e par les \'el\'ements 
de la forme $(da,0)-(0,a\otimes a)$ pour $a\in R-\{0\}$. Pour 
tout $a\in K^\times$, la classe de $(0,1\otimes a)$ sera not\'ee $d\log(a)$. 
Il est imm\'ediat de v\'erifier que l'homomorphisme
$\Omega^1_R\oplus R\rightarrow \Omega^1_R(\log)$
qui envoie $(\omega,a)$ sur $\omega+ad\log \pi$ 
est surjectif de noyau engendr\'e par $(d\pi,0)-(0,\pi)$. 
On a donc un isomorphisme canonique 
\begin{equation}\label{diff35}
\Omega^1_R\otimes_RK\stackrel{\sim}{\rightarrow} \Omega^1_R(\log)\otimes_RK.
\end{equation}

On pose $\Omega^1_F(\log)=\Omega^1_R(\log)\otimes_RF$ 
et on note $d\log[\pi]$ la classe de $d\log(\pi)$.
On voit aussit\^ot que $\Omega^1_F(\log)$ s'identifie canoniquement  
au quotient de $\Omega^1_F\oplus (F\otimes_{\mZ}K^\times)$
par le sous-$F$-module engendr\'e par les \'el\'ements 
de la forme $(d\oa,0)-(0,\oa\otimes a)$ pour $a\in R-\{0\}$ et 
$\oa$ sa classe dans $F$. On a alors une suite exacte canonique 
\begin{equation}\label{res1}
\xymatrix{
0\ar[r]&{\Omega^1_F}\ar[r]&{\Omega^1_F(\log)}
\ar[r]^-(0.5){\res}&F\ar[r]&0},
\end{equation}
o\`u $\res(0,a\otimes b)=a\ttv(b)$ pour $a\in F$ et $b\in K^\times$, 
qui est scind\'ee par l'\'el\'ement $d\log[\pi]$. 
On en d\'eduit que $\Omega^1_R(\log)$ est engendr\'e par $d$ \'el\'ements
sur $R$; il est donc libre de rang $d$ \eqref{diff35}. 

\subsection{}\label{not2}
Pour $L$ une extension finie de $K$, on pose $R_L$ 
la cl\^oture int\'egrale de $R$ dans $L$, $\eta_L=\Spec(L)$,
$S_L=\Spec(R_L)$ et $\ttv_L^{}\colon L^\times \rightarrow \mQ$ la 
valuation de $L$ qui prolonge $\ttv$ sur $K$.   
Pour $D$ un diviseur (de Cartier) effectif sur $S_L$, on pose 
$\co_D(D)=\co_{S_L}(D)\otimes_{R_L}\co_D$ et 
$\ttv_L^{}(D)=\ttv_L^{}(f)$, o\`u $f\in L$ est un g\'en\'erateur de 
$\co_{S_L}(-D)$.

\subsection{}\label{fpsi}
On d\'esigne par $\Lambda$ un anneau commutatif noeth\'erien,  
annul\'e par un entier inversible dans $k$ et tel que $\Spec(\Lambda)$
soit connexe. Soit $\psi:\mF_p\rightarrow \Lambda^\times$ 
un caract\`ere additif non-trivial. 
Si $X$ est un $k$-sch\'ema,
on note $\bD(X)=\bD(X,\Lambda)$ la cat\'egorie d\'eriv\'ee des complexes 
de $\Lambda$-modules sur $X$ pour la topologie \'etale, $\bD_\ctf(X)$ 
la sous-cat\'egorie pleine form\'ee des complexes de tor-dimension finie,
\`a cohomologie constructible. On note avec un indice $\ttc$ la 
sous-cat\'egorie pleine de $\bD(X)$ (ou $\bD^+(X)$, $\bD^-(X)$, $\bD^\ttb(X)$)
form\'ee des complexes \`a cohomologie constructible.

\subsection{}\label{CD1}
Soient $X,Y$ deux $k$-sch\'emas, $\cF$ (resp. $\cG$) un faisceau \'etale
en $\Lambda$-modules sur $X$ (resp. $Y$). On d\'esigne 
par $\cF\boxtimes\cG$ et $\cH(\cG,\cF)$ les faisceaux sur $X\times_kY$
d\'efinis par  
\begin{eqnarray*}
\cF\boxtimes\cG&=&p_1^*(\cF)\otimes_\Lambda p_2^*(\cG),\\ 
\cH(\cG,\cF)&=&\cHom(p_2^*(\cG),p_1^*(\cF)),
\end{eqnarray*}
o\`u $p_1\colon X\times_kY\rightarrow X$ et 
$p_2\colon X\times_kY\rightarrow Y$ sont les projections canoniques.
On utilisera abusivement les m\^emes notations pour designer les foncteurs 
d\'eriv\'es des foncteurs ci-dessus.

\section{Alg\'ebrisation}

\subsection{}\label{alg1}
Soient $X$ un $k$-sch\'ema affine et lisse, 
$\xi\colon S\rightarrow X$ un $k$-morphisme, $\ttx=\xi(s)$,
$A$ l'anneau local de $X$ en $\ttx$. 
On suppose que $\xi$ induit un $k$-isomorphisme 
entre les s\'epar\'es compl\'et\'es des anneaux locaux $\hA\simeq R$. 
On peut construire un tel morphisme $\xi$ de la fa\c{c}on suivante.
Comme $k$ est parfait et $F$ est de type fini sur $k$, 
il existe un $k$-sch\'ema affine et lisse  
$X_0$ de corps des fonctions rationnelles $F$; 
et comme $R$ est $k$-isomorphe \`a $F[[\pi]]$, 
on peut prendre $X=\mA^1_{k}\times_kX_0$ et $\ttx$ le point g\'en\'erique 
d'une fibre sp\'eciale de la projection canonique $X\rightarrow \mA^1_k$.  

On d\'esigne par $\fE$ la cat\'egorie des sch\'emas point\'es $(X',\ttx')$ 
au dessus de $(X,\ttx)$ tel que le morphisme structural 
$f\colon X'\rightarrow X$ soit affine, \'etale et induise 
un isomorphisme $k(\ttx')\simeq k(\ttx)$.
Pour tout objet $(X',\ttx')$ de $\fE$, 
on note $X'_0$ l'adh\'erence sch\'ematique de $\ttx'$ dans $X'$. 
On rappelle qu'il existe un et un unique $k$-morphisme 
$\xi'\colon S\rightarrow X'$ relevant $\xi$ car $S$ est hens\'elien 
(EGA IV 18.6.2). 

On note $A^\tth$ le hens\'elis\'e de $A$, 
de sorte qu'on ait un isomorphisme canonique 
\[
\Spec(A^\tth)\simeq \underset{\underset{(X',\ttx')\in \fE}
{\longleftarrow}}{\lim}\ X'.
\]
On identifie le s\'epar\'e compl\'et\'e de $A^\tth$ avec $R$ 
(EGA IV 18.6.6). On d\'esigne par $\fC(A^\tth)$ (resp. $\fC(R)$)
la cat\'egorie des $A^\tth$-alg\`ebres (resp. $R$-alg\`ebres) finies, 
plates et g\'en\'eriquement \'etales. Par compl\'etion,
o\`u ce qui revient au m\^eme par extension 
des scalaires $A^\tth\rightarrow R$, 
on obtient une \'equivalence de cat\'egories (\cite{artin} 3.11)
\begin{equation}\label{artin3.11}
\fC(A^\tth)\rightarrow \fC(R).
\end{equation}

\begin{lem}\label{alg15}
Pour tout extension finie s\'eparable $L$ de $K$,
il existe un objet $(X',\ttx')$ de $\fE$ et un morphisme fini $f\colon
Y'\rightarrow X'$ v\'erifiant les conditions suivantes~: 

{\rm (i)}\ $Y'$ est lisse sur $k$ et $f$ est \'etale au dessus de 
l'ouvert $U'$ compl\'ementaire de $X'_0$ dans $X'$;

{\rm (ii)}\ Si $\xi'\colon S\rightarrow X'$ d\'esigne le rel\`evement 
canonique de $\xi$, on a un $S$-isomorphisme  
\begin{equation}\label{is5}
Y'\times_{X'}S\simeq S_L;
\end{equation}

{\rm (iii)}\ Si de plus $L$ est une extension galoisienne finie de $K$, 
alors \eqref{is5} induit un isomorphisme entre les groupes d'automorphismes
$\Aut_{X'}(Y')\stackrel{\sim}{\rightarrow}\Aut_{K}(L)$. 
En particulier, la restriction de $f$ au dessus de $U'$ est un 
rev\^etement \'etale galoisien de groupe $\Aut_{K}(L)$. 
\end{lem}

Comme $k$ est parfait, il revient au m\^eme de demander que 
$Y'$ soit lisse ou qu'il soit r\'egulier; 
donc la condition (i) r\'esulte de (ii) quitte \`a remplacer $X'$
par un ouvert convenable contenant $\ttx'$. 
Les propositions (ii) et (iii) r\'esultent de l'\'equivalence de cat\'egories 
\eqref{artin3.11} et EGA IV 8.8.2(ii), 8.8.2.5 et 10.8.5.

\begin{lem}\label{alg2}
Pour tout $\cF\in \Ob \bD_\ctf(\eta)$, il existe un objet 
$(X',\ttx')$ de $\fE$ et si on note 
$U'$ l'ouvert compl\'ementaire de $X'_0$ dans $X'$,
$\cG\in \Ob \bD_\ctf(U')$ tel que $\xi_{U'}^{\prime *}(\cG)\simeq \cF$,
o\`u $\xi'\colon S\rightarrow X'$ est le rel\`evement canonique de $\xi$. 
\end{lem}

On rappelle que $\cF$ est isomorphe \`a un complexe born\'e de faisceaux 
constructibles en $\Lambda$-modules plats  sur $\eta$. 
Soit $L$ une extension galoisienne finie de $K$
qui trivialise chacun de ces faisceaux. Appliquant \ref{alg15}(iii),  
on en d\'eduit un complexe born\'e de faisceaux localement constants 
constructibles  en $\Lambda$-modules plats $\cG$ sur $U'$ tel que 
$\xi_{U'}^{\prime *}(\cG)\simeq \cF$.

\begin{lem}\label{local-acyc}
Le morphisme $\xi\colon S\rightarrow X$ est universellement 
localement acyclique.  
\end{lem}

Le morphisme canonique $\Spec(A^\tth)\rightarrow X$ est \'evidemment 
universellement localement acyclique. 
Le morphisme $\Spec(R)\rightarrow \Spec(A^\tth)$ 
d\'eduit de $\xi$ \'etant r\'egulier (EGA IV 6.8.1), est aussi
universellement localement acyclique (SGA 4 XIX 4.1). La preuve de 
loc. cit. d\'epend de la r\'esolution des singularit\'es en caract\'eristique
$p$, mais le th\'eor\`eme de changement de base formel de Fujiwara-Gabber 
(\cite{fujiwara} 7.1.1) permet d'\'eliminer cette condition. 
En effet, la r\'esolution des singularit\'es intervient dans la 
preuve de SGA 4 XIX 4.1 \`a travers le lemme 2.4, qu'il suffit de remplacer 
par \cite{fujiwara} 7.1.4.

\begin{lem}\label{alg3}
Soient $Y$ un sch\'ema, $f\colon Y'\rightarrow Y$ un morphisme universellement 
localement acyclique, $\cF\in \Ob \bD_\ttc^-(Y)$ et $\cG\in \Ob \bD^+(Y)$. 
Alors l'homomorphisme canonique 
\[
f^*\rR\cHom(\cF,\cG)\rightarrow \rR\cHom(f^*\cF,f^*\cG)
\]
est un isomorphisme.
\end{lem}
Il suffit de calquer SGA 5 I 4.3.

\section{Sp\'ecialisation et Micro-localisation}\label{spml}

\subsection{}\label{sml}
Dans cette section, on d\'esigne par $(X,\ttx)$ le $k$-sch\'ema point\'e 
et $\xi\colon S\rightarrow X$ le $k$-morphisme du 
\eqref{alg1}, par $L$ une extension finie s\'eparable de $K$
et par $D$ un diviseur effectif non-nul de $S_L$. 
On se donne un $K$-homomorphisme $\tau\colon L\rightarrow \oK$, 
ce qui d\'etermine deux points g\'eom\'etriques $\oeta\rightarrow S_L$ 
et $\os\rightarrow S_L$ \eqref{not1}.
On note $X_0$ l'adh\'erence sch\'ematique de $\ttx$ dans $X$ (qui est un 
diviseur de Cartier), $U$ l'ouvert compl\'ementaire de $X_0$ dans $X$, 
$f\colon S_L\rightarrow S$ le morphisme structural et
$i\colon D\rightarrow S_L$ l'injection canonique.

\subsection{}\label{sml1}
On pose $\bT_D$ le fibr\'e vectoriel $\bV(\Omega^1_R\otimes_R \co_D(D))$ 
sur $D$. Appliquons la construction \eqref{dilat16} au morphisme
$\xi\circ f\colon S_L\rightarrow X$; on d\'esigne par 
$(X\times_kS_L)_{(D)}$ la dilatation de
$X\times_kS_L$ le long du graphe de $\xi\circ f$ d'\'epaisseur $D$. 
Comme $\xi^*(\Omega^1_{X/k})\simeq \Omega^1_R$, on a 
un diagramme commutatif canonique \`a carr\'es cart\'esiens 
\begin{equation}\label{sml2} 
\xymatrix{
{\bT_D}\ar[r]\ar[d]&{(X\times_kS_L)_{(D)}}\ar[d]&
{X\times_k\eta_L}\ar[l]\ar[d]&\\
D\ar[r]&{S_L}&{\eta_L}\ar[l]}
\end{equation}

\subsection{}\label{sml3}
Appliquons la construction \eqref{dilat16} au morphisme $f$;
on d\'esigne par $(S\times_kS_L)_{(D)}$ la dilatation de $S\times_kS_L$
le long du graphe de $f$ d'\'epaisseur $D$. 
On a un diagramme commutatif canonique \`a carr\'es cart\'esiens 
\begin{equation}\label{sml4} 
\xymatrix{
{\bT_D}\ar[r]\ar[d]&{(S\times_kS_L)_{(D)}}\ar[d]&
{S\times_k\eta_L}\ar[l]\ar[d]&\\
D\ar[r]&{S_L}&{\eta_L}\ar[l]}
\end{equation}
Pour le voir, il suffit de noter que $\xi$ est plat et le diagramme 
\[
\xymatrix{
D\ar[r]^-(0.5){f\circ i}\ar@{=}[d]&S\ar[d]^\xi\\
D\ar[r]_-(0.5){\xi\circ f\circ i}&X}
\]
est cart\'esien.  Donc en vertu de \ref{dilat17},  
on a un diagramme commutatif canonique \`a carr\'es cart\'esiens 
\begin{equation}\label{sml5}
\xymatrix{
{E^S_{(D)}}\ar[r]\ar[d]_{e}&{(S\times_kS_L)_{(D)}}\ar[r]\ar[d]^{\xi_{(D)}}
&S\ar[d]^{\xi}\\
{\bT_D}\ar[r]&{(X\times_kS_L)_{(D)}}\ar[r]&X}
\end{equation}
o\`u on a not\'e $E^S_{(D)}$ le diviseur exceptionnel 
de $(S\times_kS_L)_{(D)}$; de plus $e$ est un isomorphisme. 

\begin{rema}
On peut d\'evelopper un analogue formel de \eqref{sml3}
o\`u l'on remplace $S\times_kS_L$ par son compl\'et\'e 
le long de $D$ diagonalement plong\'e.
\end{rema}

\subsection{}\label{sml6}
On pose $\bT^\ttt_D\rightarrow D$ le fibr\'e vectoriel dual de $\bT_D$, 
$\obT_D=\bT_{D}\times_D\os$, $\obT^\ttt_D=\bT_D^\ttt\times_D\os$ et 
\begin{equation}\label{fourier}
\fF\colon \bD^+(\obT_{D})\rightarrow\bD^+(\obT^\ttt_{D})
\end{equation}
la transform\'ee de Fourier-Deligne associ\'ee 
au caract\`ere $\psi$ \eqref{fpsi} (\cite{laumon} 1.2.1.1).
On appelle {\em foncteur de sp\'ecialisation le long de $D$}
et on note 
\begin{equation}\label{nu_D}
\nu_D:\bD^+(S\times_{k}\oeta)\rightarrow \bD^+(\obT_{D})
\end{equation}
le foncteur des cycles proches (\cite{for-evan} 2.1.1)
pour le morphisme $(S\times_kS_L)_{(D)}\rightarrow S_L$.
On appelle {\em foncteur de micro-localisation  le long de $D$ 
associ\'e au caract\`ere $\psi$} et on note 
\begin{equation}\label{mu_D}
\mu_D:\bD^+(S\times_{k}\oeta)\rightarrow \bD^+(\obT_{D}^\ttt)
\end{equation}
le foncteur d\'efini par $\mu_D=\fF\circ \nu_D$.

\subsection{}\label{bc}
Soient $L'$ une extension finie s\'eparable de $L$, 
$q\colon S_{L'}\rightarrow S_L$ le morphisme structural, $D'=q^{-1}(D)$. 
D'apr\`es \ref{dilat11}, on a un diagramme cart\'esien canonique 
\begin{equation}\label{bc2}
\xymatrix{
{(S\times_kS_{L'})_{(D')}}\ar[r]^{q_{(D)}}\ar[d]&
{(S\times_kS_L)_{(D)}}\ar[d]\\
{S_{L'}}\ar[r]^q&{S_L}}
\end{equation}
La restriction de $q_{(D)}$ au dessus de $\bT_D$ identifie $\bT_{D'}$ 
avec $\bT_D\times_{D}D'$. Soit $L'\rightarrow \oK$ un $K$-homomorphisme 
prolongeant $\tau$ \eqref{sml}. On en d\'eduit des isomorphismes
\[
\oq_D\colon \obT_{D'}\stackrel{\sim}{\rightarrow} \obT_{D} \ \ \
{\rm et} \ \ \
\oq^{\ttt}_{D}\colon \obT^\ttt_{D'}\stackrel{\sim}{\rightarrow} 
\obT^\ttt_{D}.
\]

\begin{prop}\label{CD2}
Soient $\cF$ un objet de $\bD_\ctf(\eta)$, 
$\cF_!$ son extension par $0$ \`a $S$.

{\rm (i)}\ Les complexes $\nu_D(\cF_!\boxtimes \Lambda)$,  
$\mu_D(\cF_!\boxtimes \Lambda)$, $\nu_D(\cH(\cF_\oeta,\cF_!))$ 
et $\mu_D(\cH(\cF_\oeta,\cF_!))$ \eqref{CD1}
sont de tor-dimension finie et leurs cohomologies sont constructibles. 

{\rm (ii)}\ Sous les hypoth\`eses de \eqref{bc}, 
on a des isomorphismes canoniques 
\begin{eqnarray*}
\oq^*_D(\nu_D(\cF_!\boxtimes\Lambda))&\stackrel{\sim}{\rightarrow}& 
\nu_{D'}(\cF_!\boxtimes\Lambda),\\
\oq^{\ttt *}_D(\mu_D(\cF_!\boxtimes\Lambda))&\stackrel{\sim}{\rightarrow}& 
\mu_{D'}(\cF_!\boxtimes\Lambda),
\end{eqnarray*}
et de m\^eme pour $\cH(\cF_\oeta,\cF_!)$. 
\end{prop}

Reprenons le diagramme \eqref{sml5} et notons 
$\Psi_D:\bD^+(X\times_k\oeta)\rightarrow \bD^+(\obT_{D})$
le foncteur des cycles proches pour le morphisme 
$(X\times_kS_L)_{(D)}\rightarrow S_L$.
Il r\'esulte de \ref{local-acyc} et SGA 4 XVI 1.1 que 
le morphisme de changement de base $\Psi_D\rightarrow \nu_D\circ \xi_{(D)}^*$
(\cite{for-evan} 2.1.7.2) est un isomorphisme. 
En vertu de \ref{alg2}, quitte \`a remplacer $(X,\ttx)$ par un objet de $\fE$, 
on peut supposer qu'il existe $\cG\in \Ob \bD_\ctf(U)$ tel que 
$\xi_U^*(\cG)\simeq \cF$. Si on note $\cG_!$ l'extension de $\cG$ 
par $0$ \`a $X$, alors 
$\cG_!\boxtimes\Lambda$ et $\cH(\cF_\oeta,\cG_!)$
appartiennent \`a $\bD_\ctf(X\times_k\oeta)$ (\cite{sga45} 1.6 et 1.7).    
D'apr\`es \ref{local-acyc} et \ref{alg3}, 
on a $(\xi\times\id_{\oeta})^*(\cH(\cF_\oeta,\cG_!))
\simeq \cH(\cF_\oeta,\cF_!)$.
La proposition (i) r\'esulte de ce qui pr\'ec\`ede, \cite{sga45} 3.2 
et \cite{for-evan} 2.1.13. La proposition (ii) r\'esulte de \cite{sga45} 3.7; 
en effet, on a un diagramme cart\'esien canonique
\[
\xymatrix{
{(X\times_kS_{L'})_{(D')}}\ar[r]\ar[d]&
{(X\times_kS_{L})_{(D)}}\ar[d]\\
{S_{L'}}\ar[r]_q&{S_L}}
\]
compatible avec \eqref{bc2}.

\begin{defi}\label{CD12}
Soient $\cF\in \Ob \bD_\ctf(\eta)$,  $\cF_!$ son extension par $0$ \`a $S$.  

{\rm (i)}\ Soient $L$ une extension finie s\'eparable de $K$, 
$D$ un diviseur effectif non-nul de $S_L$.
On pose $\ttC_{D}(\cF)$ le sous-sch\'ema r\'eduit 
de $\obT^\ttt_D$ d'espace sous-jacent le support 
de $\mu_{D}(\cH(\cF_\oeta,\cF_!))$,
c'est \`a dire la r\'eunion des supports
de ses faisceaux de cohomologie.  

{\rm (ii)}\ Soit $r>0$ un nombre rationnel. 
Soient $L$ une extension finie s\'eparable de $K$, $D$ un diviseur de $S_L$ 
de valuation $\ttv_L^{}(D)=r$ (cf. \ref{not2}). 
On pose $\obT_r=\obT_{D}$, $\obT^\ttt_r=\obT^\ttt_{D}$ et  
$\ttC_r(\cF)=\ttC_D(\cF)$.  
Ces d\'efinitions ne d\'ependent pas du choix de $L$ \`a un isomorphisme 
canonique pr\`es \eqref{CD2}.
\end{defi}

On notera que pour $\cF$ un faisceau \'etale constructible 
en $\Lambda$-modules plats sur $\eta$, 
$\ttC_{D}(\cF)$ co\"{\i}ncide avec le support 
de $\mu_{D}(\cF_!\boxtimes\Lambda)$.

\section{Sp\'ecialisation et Micro-localisation~: une variante logarithmique}
\label{lspml}

\subsection{}\label{lsml}
Conservons les hypoth\`eses de \eqref{sml}. 
On pose $\Theta_D$ le fibr\'e vectoriel 
$\bV(\Omega^1_R(\log)\otimes_R \co_D(D))$ sur $D$,
$X\times_k^{\log}S_L$ et $S\times_k^{\log}S_L$ 
les produits fibr\'es logarithmiques relativement aux diviseurs de Cartier 
$X_0$ sur $X$, $s$ sur $S$ et l'image r\'eciproque $s_L$ de $s$ sur $S_L$
\eqref{plog3}. Comme $\xi^{-1}(X_0)=s$, on peut appliquer  
la construction \eqref{plog4} au morphisme $\xi\circ f$;
on d\'esigne par $(X\times_k^{\log}S_L)_{(D)}$ la dilatation 
de $X\times_k^{\log}S_L$ le long du graphe logarithmique de 
$\xi\circ f$ d'\'epaisseur $D$. 
Le corps $k$ \'etant parfait, $X_0$ est lisse sur $k$ au voisinage de
son point g\'en\'erique $\ttx$. Donc il r\'esulte de \eqref{plog65} qu'on a 
un diagramme commutatif canonique \`a carr\'es cart\'esiens 
\begin{equation}\label{lsml2} 
\xymatrix{
{\Theta_D}\ar[r]\ar[d]&{(X\times_k^{\log}S_L)_{(D)}}\ar[d]&
{U\times_k\eta_L}\ar[l]\ar[d]&\\
D\ar[r]&{S_L}&{\eta_L}\ar[l]}
\end{equation}

\subsection{}\label{lsml3}
Appliquons la construction \eqref{plog4} au morphisme $f$;
on d\'esigne par $(S\times^{\log}_kS_L)_{(D)}$ 
la dilatation de $S\times^{\log}_kS_L$
le long du graphe logarithmique de $f$ d'\'epaisseur $D$. 
On a un diagramme commutatif canonique \`a carr\'es cart\'esiens 
\begin{equation}\label{lsml4} 
\xymatrix{
{\Theta_D}\ar[r]\ar[d]&{(S\times^{\log}_kS_L)_{(D)}}\ar[d]&
{\eta\times_k\eta_L}\ar[l]\ar[d]&\\
D\ar[r]&{S_L}&{\eta_L}\ar[l]}
\end{equation}
Pour le voir, il suffit de noter que $\xi$ est plat et le diagramme 
\[
\xymatrix{
D\ar[r]^-(0.5){f\circ i}\ar@{=}[d]&S\ar[d]^\xi\\
D\ar[r]_-(0.5){\xi\circ f\circ i}&X}
\]
est cart\'esien.  Donc en vertu de \ref{plog7},  
on a un diagramme commutatif canonique \`a carr\'es cart\'esiens 
\begin{equation}\label{lsml5}
\xymatrix{
{E^S_{(D)}}\ar[r]\ar[d]_{e}&{(S\times_k^{\log}S_L)_{(D)}}
\ar[r]\ar[d]^{\xi_{(D)}}&S\ar[d]^{\xi}\\
{\Theta_D}\ar[r]&{(X\times_k^{\log}S_L)_{(D)}}\ar[r]&X}
\end{equation}
o\`u on a not\'e $E^S_{(D)}$ le diviseur exceptionnel 
de $(S\times_k^{\log}S_L)_{(D)}$; de plus $e$ est un isomorphisme. 

\subsection{}\label{lsml6}
On pose $\Theta^\ttt_D\rightarrow D$ le fibr\'e vectoriel dual de $\Theta_D$, 
$\oTheta_D=\Theta_{D}\times_D\os$, 
$\oTheta^\ttt_D=\Theta_D^\ttt\times_D\os$ et 
\begin{equation}\label{lfourier}
\fF\colon \bD^+(\oTheta_{D})\rightarrow\bD^+(\oTheta^\ttt_{D})
\end{equation}
la transform\'ee de Fourier-Deligne associ\'ee 
au caract\`ere $\psi$ \eqref{fpsi} (\cite{laumon} 1.2.1.1).
On appelle {\em foncteur de sp\'ecialisation logarithmique
le long de $D$} et on note 
\[
\nu_D^{\log}:\bD^+(\eta\times_k\oeta)\rightarrow \bD^+(\oTheta_{D})
\]
le foncteur des cycles proches pour le morphisme $(S\times_k^{\log}S_L)_{(D)}
\rightarrow S_L$. On appelle {\em foncteur de micro-localisation logarithmique
le long de $D$ associ\'e au caract\`ere $\psi$} et on note 
\[
\mu_D^{\log}:\bD^+(\eta\times_k\oeta)\rightarrow \bD^+(\oTheta_{D}^\ttt)
\]
le foncteur d\'efini par $\mu_D^{\log}=\fF\circ \nu_D^{\log}$.

\subsection{}\label{lbc}
Soient $L'$ une extension finie s\'eparable de $L$, 
$q\colon S_{L'}\rightarrow S_L$ le morphisme structural, $D'=q^{-1}(D)$. 
D'apr\`es \ref{plog15}, on a un diagramme cart\'esien canonique 
\[
\xymatrix{
{S\times_k^{\log}S_{L'}}\ar[r]\ar[d]&
{S\times_k^{\log}S_L}\ar[d]\\
{S_{L'}}\ar[r]^q&{S_L}}
\]
et le graphe logarithmique 
de $f\circ q$ est d\'eduit de celui de $f$ par le changement de base $q$.  
Il en r\'esulte par \ref{dilat11} un diagramme cart\'esien canonique 
\begin{equation}\label{lbc2}
\xymatrix{
{(S\times_k^{\log}S_{L'})_{(D')}}\ar[r]^{q_{(D)}}\ar[d]&
{(S\times_k^{\log}S_L)_{(D)}}\ar[d]\\
{S_{L'}}\ar[r]^q&{S_L}}
\end{equation}
La restriction de $q_{(D)}$ au dessus de $\Theta_D$ identifie $\Theta_{D'}$ 
avec $\Theta_D\times_{D}D'$. Soit $L'\rightarrow \oK$ un $K$-homomorphisme 
prolongeant $\tau$ \eqref{sml}. On en d\'eduit des isomorphismes
\[
\oq_D\colon \oTheta_{D'}\stackrel{\sim}{\rightarrow} \oTheta_{D} \ \ \
{\rm et} \ \ \
\oq^{\ttt}_{D}\colon \oTheta^\ttt_{D'}\stackrel{\sim}{\rightarrow} 
\oTheta^\ttt_{D}.
\]

\begin{prop}\label{lCD2}
Soit $\cF$ un objet de $\bD_\ctf(\eta)$.

{\rm (i)}\ Les complexes $\nu^{\log}_D(\cF\boxtimes \Lambda)$,  
$\mu^{\log}_D(\cF\boxtimes \Lambda)$, $\nu^{\log}_D(\cH(\cF_\oeta,\cF))$,  
et $\mu^{\log}_D(\cH(\cF_\oeta,\cF))$ \eqref{CD1} 
sont de tor-dimension finie et leurs cohomologies sont constructibles. 

{\rm (ii)}\ Sous les hypoth\`eses de \eqref{lbc}, 
on a des isomorphismes canoniques 
\begin{eqnarray*}
\oq^*_D(\nu^{\log}_D(\cF\boxtimes\Lambda))&\stackrel{\sim}{\rightarrow}& 
\nu^{\log}_{D'}(\cF\boxtimes\Lambda),\\
\oq^{\ttt *}_D(\mu^{\log}_D(\cF\boxtimes\Lambda))
&\stackrel{\sim}{\rightarrow}& \mu^{\log}_{D'}(\cF\boxtimes\Lambda),
\end{eqnarray*}
et de m\^eme pour $\cH(\cF_\oeta,\cF)$. 
\end{prop}
Il suffit de calquer la d\'emonstration de \ref{CD2}, 
en utilisant le diagramme \eqref{lsml5} au lieu du diagramme \eqref{sml5}

\begin{defi}\label{lCD12}
Soit $\cF\in \Ob\bD_\ctf(\eta)$. 

{\rm (i)}\ Soient $L$ une extension finie s\'eparable de $K$, 
$D$ un diviseur effectif de $S_L$.
On pose $\ttC^{\log}_{D}(\cF)$ le sous-sch\'ema r\'eduit 
de $\oTheta^\ttt_D$ d'espace sous-jacent le support 
de $\mu^{\log}_{D}(\cH(\cF_\oeta,\cF))$, 
c'est \`a dire la r\'eunion des supports
de ses faisceaux de cohomologie.

{\rm (ii)}\ Soit $r>0$ un nombre rationnel. 
On fixe une extension finie s\'eparable $L$ de $K$ 
telle qu'il existe un diviseur effectif $D$ de $S_L$ 
de valuation $\ttv_L^{}(D)=r$ (cf. \ref{not2}). 
On pose $\oTheta_r=\oTheta_{D}$, $\oTheta^\ttt_r=\oTheta^\ttt_{D}$ et  
$\ttC^{\log}_r(\cF)=\ttC^{\log}_D(\cF)$.  
Ces d\'efinitions ne d\'ependent pas du choix de $L$ \`a un isomorphisme 
canonique pr\`es \eqref{lCD2}.
\end{defi}

On notera que pour $\cF$ un faisceau \'etale constructible 
en $\Lambda$-modules plats sur $\eta$, 
$\ttC^{\log}_{D}(\cF)$ co\"{\i}ncide avec le support 
de $\mu^{\log}_{D}(\cF\boxtimes\Lambda)$.

\section{Les conjectures de l'isog\'enie}\label{enonces}

\subsection{}
On d\'esigne par $(G^r)_{\mQ_{\geq 0}}$ la filtration 
de ramification sup\'erieure de $G={\rm Gal}(\oK/K)$ et par
$(G^r_{\log})_{\mQ_{\geq 0}}$ la filtration logarithmique 
d\'efinies dans \cite{as}. On renvoie \`a {\em loc. cit.} pour la d\'efinition
des sous-groupes $G^{r+}$ et $G^{r+}_{\log}$. 

\begin{defi}\label{pc}  
Soient $\cF$ un faisceau \'etale constructible
en $\Lambda$-modules sur $\eta$, $r\geq 0$ un nombre rationnel.

{\rm (i)}\ On dit que $r$ est une {\em pente critique} de $\cF$ si 
$(\cF_\oeta)^{G^r}\subsetneq (\cF_\oeta)^{G^{r+}}$. 

{\rm (ii)}\ On dit que 
$r$ est une {\em pente logarithmique critique} de $\cF$ si 
$(\cF_\oeta)^{G^r_{\log}}\subsetneq(\cF_\oeta)^{G^{r+}_{\log}}$. 
\end{defi}

On notera qu'un faisceau irr\'eductible a une unique pente critique. 

\begin{conj}\label{ciso1}
Soient $\cF$ un faisceau \'etale constructible 
en $\Lambda$-modules plats sur $\eta$, $r>1$ une pente critique de $\cF$.
On suppose que $(\cF_{\oeta})^{G^r}=0$ et $(\cF_{\oeta})^{G^{r+}}=
\cF_{\oeta}$, de sorte que $r$ soit l'unique pente critique de $\cF$. 
Alors $\ttC_r(\cF)$ est un ensemble fini de points de $\obT^\ttt_r$
ne contenant pas l'origine {\rm (cf. \ref{CD12})}.
\end{conj} 

\begin{conj}\label{ciso2}
Soient $\cF$ un faisceau \'etale constructible 
en $\Lambda$-modules plats sur $\eta$, $r>0$ une pente logarithmique critique.
On suppose que $(\cF_{\oeta})^{G^r_{\log}}=0$ et 
$(\cF_{\oeta})^{G^{r+}_{\log}}=\cF_{\oeta}$, 
de sorte que $r$ soit l'unique pente logarithmique critique de $\cF$. 
Alors $\ttC^{\log}_r(\cF)$ est un ensemble fini de points de $\oTheta^\ttt_r$
ne contenant pas l'origine {\rm (cf. \ref{lCD12})}.
\end{conj} 

Ces conjectures ne dependent pas du choix de $\psi$;  
comme $\Spec(\Lambda)$ est connexe \eqref{fpsi}, 
un autre choix aura pour effet 
de multiplier $\ttC_r(\cF)$ et  $\ttC^{\log}_r(\cF)$ 
par un entier premier \`a $p$.

\subsection{}\label{mj}
Pour tout nombre rationnel $r$, on d\'esigne par $\fm^{(r)}$ l'id\'eal 
fractionnaire de $\oR$ form\'e des $x\in \oK$ tel que $\ttv(x)\geq r$.
On pose $\obT=\bV(\Omega^1_R\otimes_R\oF)$, $\obT^\ttt$ le $\os$-fibr\'e 
vectoriel dual et 
\begin{equation}
\obT^\ttt_r\times_{\os}\bV(\fm^{(-r)}\otimes_{\oR} \oF)
\rightarrow \obT^\ttt
\end{equation}
le morphisme de multiplication.
On d\'efinit de m\^eme des variantes logarithmiques $\oTheta$, 
$\oTheta^\ttt$...

\begin{defi}\label{def-vc}
{\rm (i)}\ Sous les hypoth\`eses de \ref{ciso1}, on appelle {\em vari\'et\'e 
caract\'eristique de $\cF$} et on note $\vc(\cF)$ 
le c\^one de $\obT^\ttt$ engendr\'e par $\ttC_r(\cF)$.

{\rm (ii)}\ Sous les hypoth\`eses de \ref{ciso2}, on appelle {\em vari\'et\'e 
caract\'eristique logarithmique de $\cF$} et on note $\vc_{\log}(\cF)$ 
le c\^one de $\oTheta^\ttt$ engendr\'e par $\ttC^{\log}_r(\cF)$.
\end{defi}

\begin{rema}
La conjecture \ref{ciso1} (resp. \ref{ciso2}) implique que 
$\vc(\cF)$ (resp. $\vc_{\log}(\cF)$) est de dimension un. 
Ces c\^ones ne dependent pas du choix de $\psi$. 
\end{rema}

\subsection{}\label{H1K}
On pose 
\begin{equation}
\rH^1(K)=\underset{\underset{n\geq 1}{\longrightarrow}}{\lim}\ 
\rH^1(K,\mZ/n\mZ)\simeq \Hom(G,\mQ/\mZ)
\end{equation}
o\`u la limite inductive est prise relativement aux homomorphismes  
$\rH^1(K,\mZ/n\mZ)\rightarrow\rH^1(K,\mZ/mn\mZ)$ d\'eduits 
des homomorphismes $\cdot m\colon \mZ/n\mZ\rightarrow\mZ/mn\mZ$. 
On d\'efinit de m\^eme les sous-groupes de $\rH^1(K)$
\begin{eqnarray}
\rH^1(K)[p^\infty]&=&\underset{\underset{m\geq 1}{\longrightarrow}}{\lim}\ 
\rH^1(K,\mZ/p^{m}\mZ),\\
\rH^1(K)[p']&=&\underset{\underset{(n,p)=1}{\longrightarrow}}{\lim}\ 
\rH^1(K,\mZ/n\mZ).
\end{eqnarray}
Pour $\chi\in \rH^1(K)$, on rappelle dans \eqref{sw6} (resp. \eqref{msw6}) 
les d\'efinitions des conducteurs de Swan de Kato $\sw(\chi)$ et $\rsw(\chi)$ 
(resp. des conducteurs de Swan modifi\'es de Kato-Matsuda
$\sw'(\chi)$ et $\rsw'(\chi)$).

\subsection{}\label{ntp1}
Soient $n\geq 1$ un entier, $\chi\in \rH^1(K,\mZ/n\mZ)$. 
On suppose que $n$ est un multiple de $p$ et 
$\Lambda$ contient une racine $n$-\`eme primitive de l'unit\'e.  
On fixe un homomorphisme injectif $\mZ/n\mZ\rightarrow \Lambda^\times$ tel que 
le compos\'e $\mZ/p\mZ \rightarrow 
\mZ/n\mZ\rightarrow \Lambda^\times$ avec la multiplication par $n/p$
soit le caract\`ere $\psi$ \eqref{fpsi}. 
On note encore $\chi\colon G\rightarrow \Lambda^\times$ le caract\`ere 
induit par $\chi$ et $\cF$ le faisceau \'etale en $\Lambda$-modules 
de rang $1$ associ\'e sur $\eta$. 

\begin{teo}\label{tp1}
Les hypoth\`eses \'etant celles de \eqref{ntp1}, supposons 
de plus que $\sw'(\chi)>1$.

{\rm (i)}\ L'unique pente critique de $\cF$ est l'entier $r=\sw'(\chi)+1$. 

{\rm (ii)}\ La conjecture \eqref{ciso1} est satisfaite par $\cF$. 
Plus pr\'ecis\'ement, $\ttC_r(\cF)$ 
est le point de $\obT^\ttt_r$ d\'efini par  
\[
-\rsw'(\chi)\colon F\rightarrow \Omega^1_R\otimes_R(\fm^{-r}/\fm^{-r+1}).
\] 
\end{teo}

Ce th\'eor\`eme est d\'emontr\'e dans la section \ref{Dtp1}.

\begin{teo}\label{tp2}
Les hypoth\`eses \'etant celles de \eqref{ntp1}, supposons 
de plus que $\sw(\chi)\geq 1$.

{\rm (i)}\ L'unique pente logarithmique critique de $\cF$ est l'entier 
$r=\sw(\chi)$. 

{\rm (ii)}\ La conjecture \eqref{ciso2} est satisfaite par $\cF$. 
Plus pr\'ecis\'ement, $\ttC^{\log}_r(\cF)$ est le point de 
$\oTheta^\ttt_r$ d\'efini par 
\[
-\rsw(\chi)\colon F\rightarrow \Omega^1_R(\log)\otimes_R(\fm^{-r}/\fm^{-r+1}).
\] 
\end{teo}

Ce th\'eor\`eme est d\'emontr\'e dans la section \ref{Dtp2}. 

\begin{cor}
Soient $\chi\in \rH^1(K)\simeq\Hom(G,\mQ/\mZ)$, 
$n\geq 0$ un entier, $r$ un nombre rationnel tel que $n<r\leq n+1$. 
Alors les conditions suivantes sont \'equivalentes~:

{\rm (i)}\ $\chi\in \fil_n\rH^1(K)$, le n-\`eme cran de la filtration de Kato
\eqref{filk1}; 

{\rm (ii)}\ $\chi(G^{n+}_{\log})=0$; 

{\rm (iii)}\ $\chi(G^{r}_{\log})=0$. 
\end{cor} 

Il suffit de voir que si $\chi\not=0$,  
alors $\chi(G^{\sw(\chi)+}_{\log})=0$ et $\chi(G^{\sw(\chi)}_{\log})\not=0$.
Cela r\'esulte de \ref{tp2}(i) si $\sw(\chi)\geq 1$, et de \cite{kato1} 6.1 
et \cite{as} part I 3.15 si $\sw(\chi)=0$. 

\begin{rema}
Soient $X$ un $k$-sch\'ema lisse, $U$ un ouvert de $X$ compl\'ementaire 
d'un diviseur \`a croisements normaux stricts, 
$\cF$ un faisceau \'etale localement constant 
constructible en $\Lambda$-modules libres de rang un sur $U$.
Kato (\cite{kato2} 7.1) a montr\'e que son conducteur de Swan raffin\'e
se globalise. Le th\'eor\`eme \ref{tp2} montre alors que notre notion de 
vari\'et\'e caract\'eristique permet de reconstruire celle de Kato 
localement en les points g\'en\'eriques du diviseur de ramification 
de $\cF$. Pour une comparaison avec l'aspect global de la th\'eorie 
de Kato, on renvoie \`a \cite{ascc} Section 4. 
\end{rema}

\begin{teo}\label{parf2}
Les conjectures \eqref{ciso1} et \eqref{ciso2} sont v\'erifi\'ees si $F=k$.
\end{teo} 

Ce th\'eor\`eme est d\'emontr\'e dans  \ref{pparf2}.

\section{Ramification des caract\`eres d'Artin-Schreier-Witt}\label{thkato}

\subsection{}
Soient $m\geq 0$ un entier, $\rW_{m+1}(K)$ l'anneau des vecteurs de Witt 
de longueur $m+1$. Suivant \cite{bry,kato1}, pour tout entier $n$,  
on pose $\fil_n\rW_{m+1}(K)$ le sous-ensemble de $\rW_{m+1}(K)$ form\'e 
des \'el\'ements $(x_0,\dots,x_{m})$ tel que
\begin{equation}
p^{m-i}\ttv(x_i)\geq -n \ {\rm \text pour \ tout} \ 0\leq i\leq m.
\end{equation}
On obtient ainsi une filtration croissante exhaustive de $\rW_{m+1}(K)$
par des sous-groupes (\cite{bry} 1). On pose
\[
\Gr_n\rW_{m+1}(K)=\fil_n\rW_{m+1}(K)/\fil_{n-1}\rW_{m+1}(K).
\]  
Soit $\rV:\rW_{m+1}(K)\rightarrow \rW_{m+2}(K)$ l'homomorphisme de d\'ecalage. 
On a clairement
\begin{equation}\label{fil-dec}
\rV(\fil_n\rW_{m+1}(K))\subset \fil_n\rW_{m+2}(K).
\end{equation}

\subsection{}
Soit $\rF:\rW_{m+1}\rightarrow \rW_{m+1}$ l'endomorphisme de Frobenius. 
La suite exacte de $\eta_\et$
\[
\xymatrix{
0\ar[r]&{\mZ/p^{m+1}\mZ}\ar[r]& {\rW_{m+1}}\ar[r]^{\rF-1}& {\rW_{m+1}}\ar[r]&0}
\]
induit par passage \`a la cohomologie un homomorphisme connectant surjectif
\begin{equation}\label{deltam}
\delta_{m+1}:\rW_{m+1}(K)\rightarrow \rH^1(K,\mZ/p^{m+1}\mZ).
\end{equation}
On v\'erifie imm\'ediatement que le diagramme
\begin{equation}\label{limind1}
\xymatrix{
{\rW_{m+1}(K)}\ar[r]^-(0.5){\delta_{m+1}}\ar[d]_{\rV}&{\rH^1(K,\mZ/p^{m+1}\mZ)}
\ar[d]^{.p}\\
{\rW_{m+2}(K)}\ar[r]_-(0.5){\delta_{m+2}}&{\rH^1(K,\mZ/p^{m+2}\mZ)}}
\end{equation}
est commutatif. 

\subsection{} On pose
\begin{eqnarray}\label{fh}
\fil_n \rH^1(K,\mZ/p^{m+1}\mZ)&=&\delta_{m+1}(\fil_n \rW_{m+1}(K)).
\end{eqnarray}
On obtient ainsi une filtration croissante exhaustive
de $\rH^1(K,\mZ/p^{m+1}\mZ)$. On note
\begin{eqnarray*}
\Gr_n\rH^1(K,\mZ/p^{m+1}\mZ)&=&\fil_n \rH^1(K,\mZ/p^{m+1}\mZ)/
\fil_{n-1} \rH^1(K,\mZ/p^{m+1}\mZ).
\end{eqnarray*}

\subsection{} En vertu de \eqref{fil-dec} et \eqref{limind1}, 
\eqref{fh} d\'efinit par passage \`a la limite inductive une filtration 
croissante exhaustive
$(\fil_n\rH^1(K)[p^\infty])_{n\in \mZ}$ de $\rH^1(K)[p^\infty]$
(cf. \ref{H1K}). Pour $n\geq 0$, on pose 
\begin{equation}\label{filk1}
\fil_n\rH^1(K)=\rH^1(K)[p']+\fil_n\rH^1(K)[p^\infty].
\end{equation}
Pour $n\geq 1$, on a 
\begin{equation}
\Gr_n\rH^1(K)=\fil_n\rH^1(K)/\fil_{n-1}
\rH^1(K)\simeq \fil_n\rH^1(K)[p^\infty]/\fil_{n-1}
\rH^1(K)[p^\infty].
\end{equation}

\subsection{}
Soit $\rF:\rW_{\bullet+1}\Omega^1_K \rightarrow \rW_\bullet\Omega^1_K$
l'endomorphisme de Frobenius du complexe de de Rham-Witt. L'homomorphisme 
\begin{eqnarray}\label{Fd1}
\rF^m d: \rW_{m+1}(K)\rightarrow \Omega^1_K
\end{eqnarray}
est donn\'e par la formule 
\begin{equation}\label{Fd2}
\rF^m d(x_0,\dots,x_m)=\sum_{i=0}^m x_i^{p^{m-i}-1}dx_i.
\end{equation}
En effet, posons $x=(x_0,\dots,x_m)$ et $x_+=(x_1,\dots,x_m)$, de 
sorte que $x=x_0+\rV(x_+)$. En vertu de la relation 
$\rF d\rV=d:\rW_m(K)\rightarrow
\rW_m\Omega^1_K$  (\cite{illusie1} I (2.18.3)), on a $\rF^m d(x)= \rF^md(x_0)+
\rF^{m-1}d(x_+)$. Il suffit donc de v\'erifier que 
$\rF^md(x_0)=x_0^{p^m-1}dx_0$, ce qui r\'esulte de \cite{illusie1} I (2.18.5).

\subsection{}
On d\'efinit une filtration croissante exhaustive de $\Omega^1_K$
en posant, pour tout $n\in \mZ$, $\fil_n\Omega^1_K$ 
l'image du morphisme canonique
$\Omega^1_{R}(\log)\otimes_R \fm^{-n}\rightarrow \Omega^1_K$. 
On a alors
\begin{eqnarray}\label{grOn}
\Gr_n\Omega^1_K=\fil_n\Omega^1_K/\fil_{n-1}\Omega^1_K
\simeq \Omega^1_F(\log) \otimes_F(\fm^{-n}/\fm^{-n+1}).\nonumber
\end{eqnarray}
Il est clair que
\begin{eqnarray}\label{Fm-dec}
\rF^m d (\fil_n\rW_{m+1}(K))\subset \fil_n\Omega^1_K.
\end{eqnarray}
On en d\'eduit un homomorphisme canonique 
\begin{equation}
\gr_n(\rF^md):\Gr_n\rW_{m+1}(K)\rightarrow \Gr_n\Omega^1_K.
\end{equation}
Il r\'esulte de \eqref{fil-dec} et la relation $\rF^{m+1} d \rV=\rF^{m}d:
\rW_{m+1}(K)\rightarrow \Omega^1_K$ qu'on a 
\begin{equation}\label{limind20}
\gr_n(\rF^md) = \gr_n(\rF^{m+1}d)\circ \rV.
\end{equation}

\begin{prop}[\cite{kato1} 3.2]\label{sw1}
Pour tout entier $n\geq 1$, il existe un et un unique homomorphisme 
\begin{equation}\label{sw11}
\psi_{m,n}: \Gr_n\rH^1(K,\mZ/p^{m+1}\mZ)\rightarrow \Gr_n\Omega^1_K
\end{equation}
rendant commutatif le diagramme suivant~:
\begin{equation}
\xymatrix{
{\Gr_n\rW_{m+1}(K)}\ar[rr]^-(0.5){-\gr_n(\rF^md)}
\ar[d]_{\gr_n(\delta_{m+1})}&&
{\Gr_n\Omega^1_K}\\
{\Gr_n\rH^1(K,\mZ/p^{m+1}\mZ)}\ar[rru]_-(0.5){\psi_{m,n}}&&}
\end{equation}
\end{prop}

L'homomorphisme $\delta_{m+1}$ induit un isomorphisme 
\[
\frac{\fil_n\rW_{m+1}(K)}
{\fil_{n-1}\rW_{m+1}(K)+\fil_n\rW_{m+1}(K)\cap (\rF-1)\rW_{m+1}(K)}
\simeq \Gr_n\rH^1(K,\mZ/p^{m+1}\mZ).
\]
Il suffit donc de montrer que 
\begin{equation}\label{sw12}
\rF^md(\fil_n\rW_{m+1}(K)\cap (\rF-1)\rW_{m+1}(K))\subset 
\fil_{n-1}\Omega^1_K.
\end{equation}

Soient $s\geq 1$ un entier, $y\in \fil_{s}\rW_{m+1}(K)-\fil_{s-1}\rW_{m+1}(K)$.
On v\'erifie facilement que   
$\rF(y)\in \fil_{ps}\rW_{m+1}(K)-\fil_{ps-1}\rW_{m+1}(K)$; 
d'o\`u  $x=\rF(y)-y\in \fil_{ps}\rW_{m+1}(K)-\fil_{ps-1}\rW_{m+1}(K)$.
On conclut que si $x=\rF(y)-y\in \fil_n\rW_{m+1}(K)$ alors  
$y\in \fil_{[\frac{n}{p}]}\rW_{m+1}(K)$, 
o\`u $[-]$ est la partie enti\`ere.
La relation $d\rF=p\rF d: \rW_{m+2}(K)\rightarrow \rW_{m+1}\Omega^1_K$
(\cite{illusie1} I (2.18.2)) entra\^{\i}ne que $\rF^md(x)=-\rF^md(y)$. 
L'inclusion \eqref{sw12} r\'esulte alors de \eqref{Fm-dec} puisque
$[\frac{n}{p}]\leq n-1$ pour $n\geq 1$.

\subsection{} 
Soit $n\geq 1$ un entier. 
Il r\'esulte de \eqref{limind1}, \eqref{limind20} et \ref{sw1} 
qu'on a un triangle commutatif 
\begin{equation}\label{lim-ind50}
\xymatrix{
{\Gr_n\rH^1(K,\mZ/p^{m+1}\mZ)}\ar[d]_{.p}\ar[rr]^-(0.5){\psi_{m,n}}
&&{\Gr_n\Omega^1_K}\\
{\Gr_n\rH^1(K,\mZ/p^{m+2}\mZ)}\ar[rru]_-(0.5){\psi_{m+1,n}}&&}
\end{equation}
Le foncteur ``limite inductive'' \'etant exact \`a droite, 
on en d\'eduit un homomorphisme 
\begin{equation}\label{sw15}
\psi_n:\Gr_n\rH^1(K)\rightarrow \Gr_n\Omega^1_K.
\end{equation}

\subsection{}\label{BrZr} 
Pour tout entier $q\geq 0$, on pose $\rZ^q\Omega^\bullet_F=
\ker(d:\Omega^q_F\rightarrow \Omega^{q+1}_F)$, 
$\rB^q\Omega^\bullet_F=d(\Omega^{q-1}_F)$,
\begin{eqnarray}\label{cartier}
\rC^{-1}\colon \Omega^q_F&\stackrel{\sim}{\rightarrow}&
\rH^q(\Omega^\bullet_F)=\rZ^q\Omega^\bullet_F/\rB^q\Omega^\bullet_F
\end{eqnarray}
l'isomorphisme de Cartier inverse, $\rC\colon \rZ^q\Omega^\bullet_F
\rightarrow \Omega^q_F$ l'op\'eration de Cartier.  
Pour tout entier $r\geq 0$, on consid\`ere les sous-groupes 
\[
\rB_r\Omega^q_F\subset \rZ_r\Omega^q_F\subset \Omega^q_F
\]
d\'efinis dans \cite{illusie1} 0 (2.2.2). On note 
$\rC^r\colon \rZ_r\Omega^q_F\rightarrow \Omega^q_F$ l'op\'eration 
de Cartier it\'er\'ee $r$ fois.

\begin{lem}\label{pbase1}
Soient $(y_1,\dots,y_c)$ une $p$-base de $F$, $I_c=\{0,\dots,p-1\}^c$, 
$I'_c=I_c-\{(0,\dots,0)\}$. Pour $\uk=(k_1,\dots,k_c)\in I_c$, on pose 
$\tty^\uk=\prod_{i=1}^cy_i^{k_i}$.
Soit $r\geq 1$ un entier. Alors tout \'el\'ement
de $\rB_r\Omega^1_F$ peut s'\'ecrire d'une unique mani\`ere sous la forme 
\begin{equation}\label{pbase11}
\sum_{0\leq j\leq r-1}
\sum_{\uk\in I'_c}x_{\uk,j}^{p^{j+1}}(\tty^{\uk})^{p^j}
\sum_{1\leq i\leq c}k_i\frac{dy_i}{y_i},
\end{equation}
o\`u $x_{\uk,j}$ (pour $0\leq j\leq r-1$ et $\uk\in I_c'$) sont
des \'el\'ements de $F$. 
\end{lem}

L'\'enonc\'e est imm\'ediat pour $r=1$. 
Par d\'efinition, \eqref{cartier} induit un isomorphisme 
$\rC^{-1}\colon \rB_r\Omega^1_F\stackrel{\sim}{\rightarrow}
\rB_{r+1}\Omega^1_F/\rB_1\Omega^1_F$.
Le lemme s'en d\'eduit par r\'ecurrence sur $r$.

\subsection{}\label{Blog}
Soit $n=mp^r$ un entier $\geq 1$ tel que $r\geq 0$ et $(m,p)=1$. 
On pose $Q_n\subset \rB_{r+1}\Omega^1_F\oplus \rZ_r F$ le 
sous-groupe form\'e des \'el\'ements $(\alpha,\beta)$ tel que 
\begin{equation}\label{cond1}
m\rC^r(\alpha)=-d\rC^r(\beta).
\end{equation}
Observons qu'un \'el\'ement $\beta\in F$ appartient \`a $\rZ_r F$
si et seulement s'il existe $\gamma\in F$ tel que $\beta=\gamma^{p^r}$. 
Donc \eqref{cond1} est \'equivalente 
\`a la condition $m\rC^r(\alpha)=-d\gamma$ et on a la suite exacte
\begin{equation}\label{res2}
\xymatrix{
0\ar[r]&{\rB_r\Omega^1_F}\ar[r]^{i}&{Q_n}
\ar[r]^j&\rZ_r F\ar[r]&0}
\end{equation}
o\`u $i(\alpha)=(\alpha,0)$ et $j(\alpha,\beta)=\beta$. 

Compte tenu de \eqref{grOn}, on d\'efinit l'homomorphisme de groupes
\begin{equation}
t_{n,\pi}:Q_n\rightarrow \Omega^1_F(\log)\otimes_F
(\fm^{-n}/\fm^{-n+1})\simeq \Gr_n\Omega^1_K
\end{equation}
par la formule 
\begin{equation}
t_{n,\pi}(\alpha,\beta)=(\alpha+\beta d\log[\pi])\otimes[\pi^{-n}]
\end{equation}
o\`u $[\pi^{-n}]\in \fm^{-n}/\fm^{-n+1}$ d\'esigne la classe de $\pi^{-n}$.

\begin{lem}\label{pbase2}
Soit $n=mp^r$ un entier $\geq 1$ tel que $r\geq 0$ et $(m,p)=1$.

{\rm (i)} L'homomorphisme $t_{n,\pi}$ est injectif d'image dans 
$\Gr_n\Omega^1_K$ ind\'ependante du choix de $\pi$; 
on la note $\rB\Gr_n\Omega^1_K$. 

{\rm (ii)}\ Avec les notations de \eqref{pbase1}, 
tout \'el\'ement de $\rB\Gr_n\Omega^1_K$ 
peut s'\'ecrire d'une unique mani\`ere sous la forme 
\begin{equation}\label{pbase21}
\left(\sum_{0\leq j\leq r-1}
\sum_{\uk\in I'_c}x_{\uk,j}^{p^{j+1}}(\tty^{\uk})^{p^j}
\sum_{1\leq i\leq c}k_i\frac{dy_i}{y_i}
+x^{p^r-1}dx- mx^{p^r}d\log[\pi]\right)
\otimes [\pi^{-n}]
\end{equation}
o\`u $x_{\uk,j}$ (pour $0\leq j\leq r-1$ et $\uk\in I_c'$)
et $x$ sont des \'el\'ements de $F$.
\end{lem}
(i) L'injectivit\'e est \'evidente. 
Soient $(\alpha,\beta)\in Q_n$, $\gamma\in F$ tel que 
$\gamma^{p^r}=\beta$, $a\in F^\times$, $\ta\in R^\times$ 
un rel\`evement de $a$. On a alors 
$a^{-n}\beta =(a^{-m}\gamma )^{p^r}\in \rZ_r F$; $a^{-n}\beta d\log(a)\in 
\rZ_{r}\Omega^1_F$ (\cite{illusie1} 0 2.8.8); 
\begin{eqnarray*}
a^{-n}\alpha+a^{-n} \beta d\log[\ta\pi]&=&a^{-n}\alpha +a^{-n}\beta  d\log(a)+ 
a^{-n}\beta d\log[\pi];\\
m\rC^r(a^{-n}\alpha+ a^{-n}\beta d\log(a))&=&m a^{-m}\rC^r(\alpha)+
m a^{-m}\gamma d\log(a)\\
&=&-a^{-m}d\gamma+m\gamma a^{-m-1}da=-d(\gamma a^{-m}).
\end{eqnarray*}
On en d\'eduit que $(a^{-n}\alpha+ a^{-n}\beta d\log(a),a^{-n}\beta)
\in Q_n$ et 
\begin{equation}
t_{n,\pi}(\alpha,\beta)=t_{n,\ta\pi}(a^{-n}\alpha+ a^{-n}\beta d\log(a),
a^{-n}\beta),
\end{equation}
ce qui ach\`eve la preuve.

(ii) Cela r\'esulte de \ref{pbase1}.

\subsection{}\label{thetaM}
Soit $n=mp^r$ un entier $\geq 1$ tel que $r\geq 0$ et $(m,p)=1$. 
On pose $M_{-1}=0$ et pour tout entier $j\geq 0$, $M_{j}$ l'image de 
$\delta_{j+1}(\fil_n\rW_{j+1}(K))$ dans $\Gr_n\rH^1(K)$ \eqref{deltam}.   
Les $M_j$ forment une filtration croissante exhaustive de $\Gr_n\rH^1(K)$. 
Pour $a\in K$ et $j\geq 0$ un entier, 
on pose $\theta_j(a)=\delta_{j+1}((a,0,\dots,0))\in \rH^1(K)$.
On a clairement 
\begin{equation}\label{ap}
\theta_j(a^p)=\theta_j(a).
\end{equation}  
Pour tout $x\in M_j$, il existe 
$a\in K$ tel que $p^j\ttv(a)\geq -n$ et $x-\theta_j(a)\in M_{j-1}$. 
On en d\'eduit que $M_r=\Gr_n\rH^1(K)$, et pour tout $0\leq j\leq r$, 
$\theta_j$ induit un homomorphisme surjectif de groupes 
\begin{equation}\label{thetaM1}
\fm^{-np^{-j}}/\fm^{-np^{-j}+1}\rightarrow M_j/M_{j-1}.
\end{equation}

Si $r\geq 1$, $0\leq j\leq r-1$ et $a\in \fm^{-np^{-j-1}}$, alors  
$\theta_j(a^p)=\theta_j(a)\in \fil_{n/p}\rH^1(K)\subset \fil_{n-1}\rH^1(K)$;
d'o\`u la classe de $\theta_j(a^p)$ est nulle dans $M_j$.

\begin{prop}[\cite{kato1} 3.2 et 3.3] \label{sw3}
Pour tout entier $n\geq 1$,  
le morphisme $\psi_n$ \eqref{sw15} induit un isomorphisme 
\[
\psi_n\colon \Gr_n\rH^1(K)\stackrel{\sim}{\rightarrow}\rB\Gr_n\Omega^1_K.
\]
\end{prop}

On \'ecrit $n=mp^r$ avec $r\geq 0$ et $(m,p)=1$. 
On fixe pour tout \'el\'ement $x\in F$ un rel\`evement $\tx\in R$. 
Reprenons les notations de \eqref{pbase1}, 
de plus, pour $\uk=(k_1,\dots,k_c)\in I_c$, 
posons $\ttty^\uk=\prod_{i=1}^c\ty_i^{k_i}$.
En vertu de \ref{pbase2}, il existe une et une unique application 
\[
\rho_n:\rB\Gr_n\Omega^1_K\rightarrow \Gr_n\rH^1(K)
\]
qui envoie l'\'el\'ement \eqref{pbase21} sur la classe de (cf. \ref{thetaM})
\begin{equation}\label{pbase22}
\sum_{0\leq j\leq r-1}\sum_{\uk\in I'_c}
\theta_j(\tx^p_{\uk,j} \ttty^\uk\pi^{-np^{-j}})+\theta_r(\tx\pi^{-m}).
\end{equation}
Il r\'esulte de \ref{thetaM} que l'application $\rho_n$ est surjective;
noter que dans \eqref{pbase22}, 
on peut se limiter \`a des sommes pour $\uk\in I'_c$ (plut\^ot
que $\uk\in I_c$). 

Pour $0\leq j\leq r-1$, $a\in R$ et $\oa$ la classe de $a$ dans $F$, 
on a les relations
\begin{eqnarray}
\gr_{n}(\rF^jd)((a\pi^{-np^{-j}},0,\dots,0))&=&(\oa^{p^j-1}d \oa)
\otimes[\pi^{-n}]\\
\gr_{n}(\rF^rd)((a\pi^{-m},0,\dots,0))&=&(\oa^{p^r-1}d\oa-m\oa^{p^r}d\log[\pi])
\otimes[\pi^{-n}]
\end{eqnarray}
dans $\rB\Gr_n\Omega^1_K$. On en d\'eduit aussit\^ot que 
$\psi_n\circ \rho_n$ est l'identit\'e de $\rB\Gr_n\Omega^1_K$.
Comme $\rho_n$ est surjectif, $\psi_n$ et $\rho_n$ sont  
des isomorphismes inverses l'un de l'autre. En particulier, $\rho_n$
est un homomorphisme de groupes qui ne d\'epend d'aucun choix. 

\begin{rema}\label{sw35}
On n'utilisera de \ref{sw3} que l'injectivit\'e de $\psi_n$. 
\end{rema}

\begin{defi}\label{sw6}
Soit $\chi\in \rH^1(K)$. On appelle {\em conducteur de Swan de $\chi$}
et on note $\sw(\chi)$ le plus petit entier $n\geq 0$ tel que 
$\chi\in \fil_n\rH^1(K)$. 

Si $\sw(\chi)\geq 1$, on appelle 
{\em conducteur de Swan raffin\'e de $\chi$} et on note $\rsw(\chi)$
l'image de la classe de $\chi$ par l'homomorphisme
\[
\psi_{\sw(\chi)}\colon \Gr_{\sw(\chi)}\rH^1(K)\rightarrow \Gr_{\sw(\chi)}
\Omega^1_K.  
\]
\end{defi}

\section{Ramification des caract\`eres d'Artin-Schreier-Witt~: 
la variante de Matsuda}\label{thkato-mats}

\subsection{} Soit $m\geq 0$ un entier. 
Suivant \cite{mats}, pour tout entier $n\geq 0$,  on pose
\begin{equation}
\fil'_n\rW_{m+1}(K)=
\rV^{m+1-m'}(\fil_{n+1}\rW_{m'}(K))+\fil_n\rW_{m+1}(K),
\end{equation}
o\`u $m'=\min(\ord_p(n+1),m+1)$. Il r\'esulte de \eqref{fil-dec}
que $(\fil'_n\rW_{m+1}(K))_{n\in \mN}$ est une filtration croissante exhaustive
de $\rW_{m+1}(K)$.  Pour $n\geq 1$, on note 
\[
\Gr'_n\rW_{m+1}(K)=\fil'_n\rW_{m+1}(K)/\fil'_{n-1}\rW_{m+1}(K).
\]

\begin{lem}\label{fil'-dec}
Pour tout entier $n\geq 0$, on a $\rV(\fil'_n\rW_{m+1}(K))\subset 
\fil'_n\rW_{m+2}(K)$.
\end{lem}

Si $\min(\ord_p(n+1),m+1)=\min(\ord_p(n+1),m+2)$, 
l'assertion r\'esulte aussit\^ot de \eqref{fil-dec}. 
Sinon, $\ord_p(n+1)>m+1$, auquel cas 
$\fil'_n\rW_{m+1}(K)=\fil_{n+1}\rW_{m+1}(K)$ et $\fil'_n\rW_{m+2}(K)
=\fil_{n+1}\rW_{m+2}(K)$; donc l'assertion r\'esulte encore 
de \eqref{fil-dec}.

\subsection{} Pour tout entier $n\geq 0$, on pose (cf. \eqref{deltam})
\begin{equation}\label{fh'}
\fil'_n \rH^1(K,\mZ/p^{m+1}\mZ)=\delta_{m+1}(\fil'_n \rW_{m+1}(K)).
\end{equation}
On obtient ainsi une filtration croissante exhaustive
de $\rH^1(K,\mZ/p^{m+1}\mZ)$. Pour $n\geq 1$, on note
\[
\Gr'_n\rH^1(K,\mZ/p^{m+1}\mZ)=\fil'_n \rH^1(K,\mZ/p^{m+1}\mZ)/
\fil'_{n-1} \rH^1(K,\mZ/p^{m+1}\mZ).
\]

\subsection{}
En vertu de \eqref{limind1} et \ref{fil'-dec}, 
\eqref{fh'} d\'efinit par passage \`a la limite inductive une filtration 
croissante exhaustive 
$(\fil'_n\rH^1(K)[p^\infty])_{n\in \mN}$ de $\rH^1(K)[p^\infty]$
(cf. \ref{H1K}). Pour $n\geq 0$, on pose 
\begin{equation}\label{filkm1}
\fil'_n\rH^1(K)=\rH^1(K)[p']+\fil'_n\rH^1(K)[p^\infty].
\end{equation}
Pour $n\geq 1$, on a 
\begin{equation}
\Gr'_n\rH^1(K)=\fil'_n\rH^1(K)/\fil'_{n-1}
\rH^1(K)\simeq \fil'_n\rH^1(K)[p^\infty]/\fil'_{n-1}
\rH^1(K)[p^\infty].
\end{equation} 

\subsection{}\label{diff2}
On d\'efinit une filtration croissante exhaustive de $\Omega^1_K$
en posant, pour tout $n\in \mZ$, 
$\fil'_n\Omega^1_K$ l'image du morphisme canonique
$\Omega^1_{R}\otimes_R \fm^{-n-1}\rightarrow \Omega^1_K$. 
On a alors  
\begin{equation}\label{Grmn}
\Gr'_n\Omega^1_K=\fil'_n\Omega^1_K/\fil'_{n-1}\Omega^1_K
\simeq \Omega^1_R \otimes_R(\fm^{-n-1}/\fm^{-n}).
\end{equation}

\begin{lem}\label{lw1}
Pour tout entier $n\geq 0$, on a  
\[
\rF^m d (\fil'_n\rW_{m+1}(K))\subset \fil'_n\Omega^1_K.
\]  
\end{lem}

Comme $d\log(x)\in \fm^{-1}\Omega^1_R$ pour tout $x\in K^\times$, 
la relation \eqref{Fd2} entra\^{\i}ne que 
\begin{eqnarray}\label{lw11}
\rF^m d (\fil_n\rW_{m+1}(K))\subset \fil'_n\Omega^1_K.
\end{eqnarray}
Posons $m'=\min(\ord_p(n+1),m+1)$ et montrons 
que $\rF^m d \rV^{m+1-m'}(\fil_{n+1}\rW_{m'}(K))\subset \fil'_n\Omega^1_K$. 
On peut se borner au cas o\`u $m'\geq 1$. Compte tenu de la relation 
$\rF^m d \rV^{m+1-m'}=\rF^{m'-1}d:\rW_{m'}(K)\rightarrow \Omega^1_K$,
il suffit de montrer que si $m+1\leq \ord_p(n+1)$, alors 
\begin{eqnarray}\label{lw12}
\rF^m d (\fil_{n+1}\rW_{m+1}(K))\subset \fil'_n\Omega^1_K.
\end{eqnarray}
Soit $(x_0,\dots,x_{m})\in \fil_{n+1}\rW_{m+1}(K)$. 
Si $p^{m-i}\ttv(x_i)=-n-1$ pour un entier $0\leq i\leq m$, alors  
$\ttv(x_i)$ est un multiple de $p$ car $\ord_p(n+1)\geq m+1$; 
donc  $d\log(x_i)=0$. L'inclusion \eqref{lw12} s'ensuit 
compte tenu de \eqref{lw11}.

\subsection{}
En vertu de \ref{lw1}, pour tout $n\geq 1$, $\rF^md$ induit 
un homomorphisme canonique 
\begin{equation}
\gr'_n(\rF^md):\Gr'_n\rW_{m+1}(K)\rightarrow \Gr'_n\Omega^1_K.
\end{equation}
Il r\'esulte de \ref{fil'-dec} et la relation $\rF^{m+1} d \rV=\rF^{m}d:
\rW_{m+1}(K)\rightarrow \Omega^1_K$ qu'on a 
\begin{equation}\label{limind2}
\gr'_n(\rF^md) = \gr'_n(\rF^{m+1}d)\circ \rV.
\end{equation}

\begin{prop}[\cite{mats} 3.2.1] \label{msw1}
Pour tout entier $n>1$ (resp. $n\geq 1$ si $p\not=2$), 
il existe un et un unique homomorphisme 
\begin{equation}\label{msw11}
\phi_{m,n}: \Gr'_n\rH^1(K,\mZ/p^{m+1}\mZ)\rightarrow \Gr'_n\Omega^1_K
\end{equation}
rendant commutatif le diagramme suivant~:
\begin{equation}
\xymatrix{
{\Gr'_n\rW_{m+1}(K)}\ar[rr]^-(0.5){-\gr'_n(\rF^md)}
\ar[d]_{\gr'_n(\delta_{m+1})}&&
{\Gr'_n\Omega^1_K}\\
{\Gr'_n\rH^1(K,\mZ/p^{m+1}\mZ)}\ar[rru]_-(0.5){\phi_{m,n}}&&}
\end{equation}
\end{prop}

L'homomorphisme $\delta_{m+1}$ \eqref{deltam} induit un isomorphisme 
\[
\frac{\fil'_n\rW_{m+1}(K)}
{\fil'_{n-1}\rW_{m+1}(K)+\fil'_n\rW_{m+1}(K)\cap (\rF-1)\rW_{m+1}(K)}
\simeq \Gr'_n\rH^1(K,\mZ/p^{m+1}\mZ).
\]
Il suffit donc de montrer que 
\begin{equation}\label{msw12}
\rF^md(\fil'_n\rW_{m+1}(K)\cap (\rF-1)\rW_{m+1}(K))\subset 
\fil'_{n-1}\Omega^1_K.
\end{equation}
Soit $y\in \rW_{m+1}(K)$ tel que 
$x=\rF(y)-y\in \fil'_n\rW_{m+1}(K)$. Comme $x\in \fil_{n+1}\rW_{m+1}(K)$, 
alors $y\in \fil_{[\frac{n+1}{p}]}\rW_{m+1}(K)\subset 
\fil'_{[\frac{n+1}{p}]}\rW_{m+1}(K)$, o\`u $[-]$ est la partie enti\`ere
(voir preuve de \ref{sw1}).
La relation $d\rF=p\rF d: \rW_{m+2}(K)\rightarrow \rW_{m+1}\Omega^1_K$
(\cite{illusie1} I (2.18.2)) entra\^{\i}ne que $\rF^md(x)=-\rF^md(y)$. 
L'inclusion \eqref{msw12} r\'esulte alors de \ref{lw1} puisque
$[\frac{n+1}{p}]\leq n-1$ sous les hypoth\`eses de la proposition.

\subsection{} 
Soit $n>1$ un entier  (resp. $n\geq 1$ si $p\not=2$). 
Il r\'esulte de \eqref{limind1}, \eqref{limind2} et \ref{msw1} 
qu'on a un triangle commutatif 
\begin{equation}\label{lim-ind5}
\xymatrix{
{\Gr'_n\rH^1(K,\mZ/p^{m+1}\mZ)}\ar[d]_{.p}\ar[rr]^-(0.5){\phi_{m,n}}
&&{\Gr'_n\Omega^1_K}\\
{\Gr'_n\rH^1(K,\mZ/p^{m+2}\mZ)}\ar[rru]_-(0.5){\phi_{m+1,n}}&&}
\end{equation}
Le foncteur ``limite inductive'' \'etant exact \`a droite, 
on en d\'eduit un homomorphisme 
\begin{equation}\label{msw15}
\phi_n:\Gr'_n\rH^1(K)\rightarrow \Gr'_n\Omega^1_K.
\end{equation}

\subsection{}
Soient $n\geq 1$ un entier, $r=\ord_p(n+1)$, $r'=\ord_p(n)$.
On pose $\rB\Gr'_n\Omega^1_K$ le sous-groupe de $\Gr'_n\Omega^1_K\simeq
\Omega^1_R\otimes_R\fm^{-n-1}/\fm^{-n}$ \eqref{Grmn}
engendr\'e par les \'el\'ements de l'une des formes suivantes 
\begin{eqnarray*}
a^{p^j-1}da\otimes [\pi^{-n-1}] \  \ ({\rm pour}\ 0\leq j\leq r-1),
&&a^{p^{r'}}d\pi\otimes [\pi^{-n-1}]
\end{eqnarray*}
o\`u $a\in R$ et $[\pi^{-n-1}]$ d\'esigne la classe de $\pi^{-n-1}$. 
Pour $0\leq j\leq r-1$, $a\in R$ et $b\in R^\times$, on a 
\begin{eqnarray*}
b^{n+1}a^{p^j-1}da&=&(b^{(n+1)p^{-j}}a)^{p^j-1}d(b^{(n+1)p^{-j}}a)\\
b^{n+1}a^{p^{r'}}d(b^{-1}\pi)&=&(b^{np^{-r'}}a)^{p^{r'}}d\pi \ \mod \pi.
\end{eqnarray*}
Donc $\rB\Gr'_n\Omega^1_K$ est ind\'ependant du choix de $\pi$. 

\begin{lem}\label{pbase3}
Soient $n\geq 1$ un entier, $r=\ord_p(n+1)$, $r'=\ord_p(n)$.
Soient $(b_1,\dots,b_c)$ des \'el\'ements de $R$ 
relevant une $p$-base de $F$, $I_c=\{0,\dots,p-1\}^c$, 
$I'_c=I_c-\{(0,\dots,0)\}$. Pour $\uk=(k_1,\dots,k_c)\in I_c$, on pose 
$\ttb^\uk=\prod_{i=1}^cb_i^{k_i}$.
On fixe pour chaque 
\'el\'ement $x$ de $F$ un rel\`evement $\tx$ dans $R$.
Alors tout \'el\'ement de $\rB\Gr'_n\Omega^1_K$ 
peut s'\'ecrire d'une unique mani\`ere sous la forme 
\begin{equation}\label{pbase31}
\left(\sum_{0\leq j\leq r-1}
\sum_{\uk\in I'_c}\tx_{\uk,j}^{p^{j+1}}(\ttb^{\uk})^{p^j}
\sum_{1\leq i\leq c}k_i\frac{db_i}{b_i}
+\tx^{p^{r'}}d\pi\right)
\otimes [\pi^{-n-1}]
\end{equation}
o\`u $x_{\uk,j}$ (pour $0\leq j\leq r-1$ et $\uk\in I_c'$)
et $x$ sont des \'el\'ements de $F$. 
\end{lem}

L'assertion est imm\'ediate si $r=0$. Supposons $r\geq 1$; donc $r'=0$ et 
on a une suite exacte (d\'ependante du choix de $\pi$)
\begin{equation}
0\rightarrow F\rightarrow \rB\Gr'_n\Omega^1_K \rightarrow
\rB_r\Omega^1_F \rightarrow 0.
\end{equation}
La proposition r\'esulte alors de \ref{pbase1}. 

\subsection{}\label{thetaM'}
Soient $n\geq 1$ un entier, $r=\ord_p(n+1)$, $r'=\ord_p(n)$. On a 
\begin{eqnarray*}
\fil'_{n}\rW_{j+1}(K)&=&\fil_{n+1}\rW_{j+1}(K) \ \ 
({\rm pour}\ 0\leq j\leq r-1),\\
\fil'_{n}\rW_{j+1}(K)&=&\fil_n\rW_{j+1}(K)+\rV^{j+1-r}\fil_{n+1}\rW_{r}(K)
\ \ ({\rm pour}\ j\geq r).
\end{eqnarray*}
On pose $M'_{-1}=0$ et pour tout entier $j\geq 0$, $M'_{j}$ l'image de 
$\delta_{j+1}(\fil'_n\rW_{j+1}(K))$ dans $\Gr'_n\rH^1(K)$ \eqref{deltam}.   
Les $M'_j$ forment une filtration croissante exhaustive de $\Gr'_n\rH^1(K)$
\eqref{fil'-dec}. On rappelle que pour $a\in K$ et $j\geq 0$, on a not\'e 
$\theta_j(a)=\delta_{j+1}((a,0,\dots,0))\in \rH^1(K)$ \eqref{thetaM}.

Supposons d'abord que $r\geq 1$; donc $r'=0$ et $M'_{r-1}=\Gr'_n\rH^1(K)$.
Si $0\leq j\leq r-1$ et $a\in \fm^{-(n+1)p^{-j}}$, 
alors $\theta_j(a)\in \fil'_n\rH^1(K)$ et sa classe 
dans $\Gr'_n\rH^1(K)$ appartient \'evidemment \`a $M'_j$. 
On voit aussit\^ot que $\theta_j$ induit un homomorphisme surjectif de groupes 
\begin{equation}\label{thetaM'1}
\fm^{-(n+1)p^{-j}}\rightarrow M'_j/M'_{j-1}.
\end{equation}
On a $\fil'_{n-1}\rW_{j+1}(K)=\fil_{n-1}\rW_{j+1}(K)$ pour tout $j\geq 0$. 
Donc le noyau de \eqref{thetaM'1} contient $\fm^{-(n+1)p^{-j}+1}$ 
si $1\leq j\leq r-1$, et $\fm^{-n+1}$ si $j=0$.

Supposons ensuite que $r\geq 1$ et ($n>1$ pour $p=2$); 
donc $n+1\leq p(n-1)$. Si $0\leq j\leq r-1$ et $a\in \fm^{-(n+1)p^{-j-1}}$, 
alors $\theta_j(a^p)=\theta_j(a)$ \eqref{ap} et $\theta_j(a)\in 
\fil'_{n-1}\rH^1(K)$; d'o\`u la classe de $\theta_j(a^p)$ 
est nulle dans $M'_j$.

Supposons enfin que $r=0$; alors $M'_{r'}=\Gr'_n\rH^1(K)$. On a 
\begin{eqnarray*}
\fil'_{n-1}\rW_{j+1}(K)&=&\fil_{n}\rW_{j+1}(K)=\fil'_{n}\rW_{j+1}(K)
\ \ ({\rm pour}\ 0\leq j\leq r'-1),\\
\fil'_{n-1}\rW_{r'+1}(K)&=&\rV(\fil_n\rW_{r'}(K))+\fil_{n-1}\rW_{r'+1}(K).
\end{eqnarray*}
On en d\'eduit que $M'_j=0$ pour $0\leq j\leq r'-1$, et 
l'application d\'eduite de $\theta_{r'}$
\begin{equation}\label{thetaM'2}
\fm^{-np^{-r'}}/\fm^{-np^{-r'}+1}\rightarrow M'_{r'}
\end{equation}
est un homomorphisme surjectif de groupes. 

\begin{prop}[\cite{mats} 3.2.3] \label{msw3}
Pour tout entier $n>1$ (resp. $n\geq 1$ si $p\not=2$), 
le morphisme $\phi_n$ \eqref{msw15} induit un isomorphisme 
\[
\phi_n\colon \Gr'_n\rH^1(K)\stackrel{\sim}{\rightarrow}\rB\Gr'_n\Omega^1_K.
\]
\end{prop}

On pose $r=\ord_p(n+1)$ et $r'=\ord_p(n)$.   
Avec les notations de \ref{pbase3}, il existe une et une unique application 
\[
\kappa_n\colon \rB\Gr'_n\Omega^1_K\rightarrow \Gr_n\rH^1(K)
\]
qui envoie l'\'el\'ement \eqref{pbase31} sur la classe de (cf. \ref{thetaM'})
\begin{equation}\label{pbase32}
\sum_{0\leq j\leq r-1}\sum_{\uk\in I'_c}
\theta_j(\tx^p_{\uk,j} \ttb^\uk\pi^{-(n+1)p^{-j}})-\frac{1}{np^{-r'}}
\theta_{r'}(\tx\pi^{-np^{-r'}}).
\end{equation}
Il r\'esulte de \ref{thetaM'} que l'application $\kappa_n$ est surjective;
on distinguera les cas $r=0$, et $r\geq 1$  qui entra\^{\i}ne que $r'=0$.

Pour $0\leq j\leq r-1$ et $a\in R$, 
on a les relations
\begin{eqnarray}
\gr'_{n}(\rF^jd)((a\pi^{-(n+1)p^{-j}},0,\dots,0))&=&(a^{p^j-1}da)
\otimes[\pi^{-n-1}]\\
\gr'_{n}(\rF^{r'}d)((a\pi^{-np^{-r'}},0,\dots,0))&=&
-(np^{-r'}a^{p^{r'}}d\pi)\otimes[\pi^{-n-1}]
\end{eqnarray}
dans $\rB\Gr_n\Omega^1_K$. On en d\'eduit aussit\^ot que 
$\phi_n\circ \kappa_n$ est l'identit\'e de $\rB\Gr'_n\Omega^1_K$.
Comme $\kappa_n$ est surjectif, $\phi_n$ et $\kappa_n$ sont  
des isomorphismes inverses l'un de l'autre. En particulier, $\kappa_n$
est un homomorphisme de groupes qui ne d\'epend d'aucun choix.

\begin{rema}\label{msw35}
On n'utilisera de \ref{msw3} que l'injectivit\'e de $\phi_n$. 
\end{rema}

\begin{defi}\label{msw6}
Soit $\chi\in \rH^1(K)$. On appelle {\em conducteur de Swan modifi\'e 
de $\chi$} et on note $\sw'(\chi)$ le plus petit entier $n\geq 0$ tel que 
$\chi\in \fil'_n\rH^1(K)$. 

Si $\sw(\chi)>1$ (resp. $\sw(\chi)\geq 1$ et $p\not=2$), on appelle 
{\em conducteur de Swan modifi\'e raffin\'e de $\chi$} et on note $\rsw'(\chi)$
l'image de la classe de $\chi$ par l'homomorphisme
\[
\phi_{\sw'(\chi)}\colon \Gr'_{\sw'(\chi)}\rH^1(K)\rightarrow \Gr'_{\sw'(\chi)}
\Omega^1_K.  
\]
\end{defi}

\section{D\'emonstration du th\'eor\`eme \ref{tp1}}\label{Dtp1}

\subsection{}
Nous d\'emontrons en premier lieu des r\'esultats pr\'eliminaires sur 
les vecteurs de Witt. On d\'efinit une suite de polyn\^omes 
$Q_n\in \mZ[\frac{1}{p}][X_1,\dots,X_i,Y_0,\dots,Y_i]$ (pour $n\in \mN$) par 
la relation de r\'ecurrence   
\begin{equation}
\sum_{i=0}^np^i(X_i(1+Y_i))^{p^{n-i}}=\sum_{i=0}^np^iX_i^{p^{n-i}}
+\sum_{i=0}^np^iQ_i^{p^{n-i}}.
\end{equation}
Le lemme suivant est imm\'ediat. 

\begin{lem}\label{witt} 
{\rm (i)}\  Le polyn\^ome $Q_n$ appartient \`a l'id\'eal de
$\mZ[X_1,\dots,X_i,Y_0,\dots,Y_i]$
engendr\'e par $(Y_0,\dots,Y_n)$, et la diff\'erence 
$Q_n-\sum_{i=0}^nX_i^{p^{n-i}}Y_i$ 
appartient \`a l'id\'eal engendr\'e par $(Y_iY_j)_{0\leq i,j\leq n}$. 

{\rm (ii)}\ Si on affecte \`a la variable $X_i$ le poids $p^i$ et 
\`a la variable $Y_i$ le poids $0$, le polyn\^ome $Q_n$ devient homog\`ene 
de degr\'e $p^n$. 

{\rm (iii)}\ Soient $A$ un anneau commutatif, $x=(x_0,\dots,x_n)$,
$y=(y_0,\dots,y_n)$ deux \'el\'ements de $A^{n+1}$. Posons $x'_i=x_i(1+y_i)$
et $x'=(x'_0,\dots,x'_n)$. On a alors la relation suivante dans $\rW_{n+1}(A)$
\[
x'-x=(Q_0(x,y),Q_1(x,y),\dots,Q_n(x,y)).
\]
\end{lem}

\subsection{}\label{ig1}
Dans cette section, on d\'esigne par $(X,\ttx)$ le $k$-sch\'ema point\'e 
et $\xi\colon S\rightarrow X$ le $k$-morphisme 
du \eqref{alg1}, par $D$ un diviseur effectif de $S$. 
On note $X_0$ l'adh\'erence sch\'ematique de $\ttx$ dans $X$
(qui est un diviseur de Cartier) et 
$U$ l'ouvert compl\'ementaire de $X_0$ dans $X$. 
On pose $r=\ttv(D)$ \eqref{not2},  
$\bT_D=\bV(\Omega^1_R\otimes_R\co_D(D))$ et $\obT_D=\bT_D\times_D\os$. 
Appliquons la construction \eqref{dilat16} au morphisme
$\xi$; on d\'esigne par $(X\times_kS)_{(D)}$ la dilatation de
$X\times_kS$ le long du graphe de $\xi$ d'\'epaisseur $D$. 
On a alors un diagramme commutatif canonique \`a carr\'es cart\'esiens 
\begin{equation}\label{ig11} 
\xymatrix{
{\bT_D}\ar[r]\ar[d]&{(X\times_kS)_{(D)}}\ar[d]&
{X\times_k\eta}\ar[l]\ar[d]&\\
D\ar[r]&{S}&{\eta}\ar[l]}
\end{equation} 
On note $\kappa$ le point g\'en\'erique de $\bT_D$, 
$R_\mK$ le s\'epar\'e compl\'et\'e de l'anneau local 
de $(X\times_kS)_{(D)}$ en $\kappa$ 
(qui est un anneau de valuation discr\`ete), 
$\mK$ son corps des fractions, 
$\fm_\mK$ son id\'eal maximal. La projection canonique 
$(X\times_kS)_{(D)}\rightarrow S$ \'etant lisse,  
elle induit un homomorphisme $u\colon R\rightarrow R_\mK$
v\'erifiant $\fm_\mK=u(\fm)R_\mK$. On note encore
$\ttv\colon \mK^\times \rightarrow \mZ$ la valuation qui prolonge 
$\ttv$ sur $K$ \eqref{not1}. Compte tenu de l'interpr\'etation modulaire
de $\bT_D$, le point $\kappa$ est uniquement d\'etermin\'e 
par un $F$-homomorphisme 
\begin{equation}\label{der-t1}
t\colon \Omega^1_R\rightarrow \fm_\mK^{r}/\fm_\mK^{r+1}.
\end{equation}
Appliquant le produit tensoriel $\otimes_F(\fm^{-n-1}/\fm^{-n})$ 
(pour $n\in \mZ$), on obtient un $F$-homomorphisme
\begin{equation}
t_n\colon \Gr'_n\Omega^1_K\rightarrow \fm_\mK^{-n+r-1}/\fm_\mK^{-n+r}
=\Gr_{n-r+1}\rW_1(\mK).
\end{equation}
La premi\`ere projection $(X\times_kS)_{(D)}\rightarrow X$ envoie $\kappa$
sur $\ttx$. On en d\'eduit un second homomorphisme 
$v\colon R\rightarrow R_\mK$.

\begin{lem}\label{vx}
{\rm (i)}\ Pour tout $x\in R$, on a $v(x)-u(x)\in \fm_\mK^r$ et
\begin{equation}\label{vx1}
v(x)-u(x)\equiv t(d(x))\mod \fm_\mK^{r+1}.
\end{equation}  

{\rm (ii)}\ Pour tout $x\in K^\times$, on a  
\begin{equation}\label{vx2}
\ttv(\frac{v(x)}{u(x)}-1)\geq r-1.
\end{equation}  

{\rm (iii)}\ Supposons $r\geq 2$. Alors on a, pour tout $x\in K^\times$,
\begin{equation}\label{vx3}
\frac{v(x)}{u(x)}-1\equiv t_0(d\log(x))\mod \fm_\mK^{r}.
\end{equation}
\end{lem}
L'assertion (i) est \'evidente. L'assertion (ii) est \'evidente pour 
les unit\'es et les uniformisantes de $R$; elle r\'esulte en g\'en\'eral 
de la relation, pour tout $x,y\in K$,   
\begin{equation}\label{vx4}
v(xy)-u(xy)=v(x)(v(y)-u(y))+u(y)(v(x)-u(x)).
\end{equation}
Il r\'esulte encore de cette relation que si \eqref{vx3}
est satisfaite par $x$ et $y$, elle l'est aussi par $1/x$ et $xy$,
car $r\geq 2$. Comme \eqref{vx3} est satisfaite par 
les unit\'es et les uniformisantes de $R$, elle est alors
satisfaite en g\'en\'eral.

\begin{lem} \label{ver}
Supposons $r\geq 3$. 
Soient $m,n\geq 0$ deux entiers, $x=(x_0,\dots,x_m)\in \fil'_{n}\rW_{m+1}(K)$, 
$y=(y_0,\dots,y_m)\in \mK^{m+1}$, o\`u 
$y_i=\frac{v(x_i)}{u(x_i)}-1$ si $ x_i\not=0$ et 
$y_i=0$ si $x_i=0$. Alors on a
\begin{eqnarray}
p^{m-i}\ttv(Q_i(u(x),y))&\geq& -n+r, \ \ \ \ \ \ \ \forall 
\ 0\leq i<m,\label{ver11}\\
\ttv(Q_m(u(x),y))&\geq &-n+r-1,\label{ver21}\\
\ttv(Q_m(u(x),y)- \sum_{j=0}^mx_j^{p^{m-j}}y_j)&\geq& -n+r.\label{ver31}
\end{eqnarray}
\end{lem}
Soit $0\leq i\leq m$ un entier. Si $x\in \fil_{n}\rW_{m+1}(K)$, alors
\[
p^{m-i}\ttv(Q_i(u(x),y))\geq -n+p^{m-i}(r-1)
\]
en vertu de \eqref{vx2} et \ref{witt}, ce qui d\'emontre \eqref{ver11}
et \eqref{ver21}.  De m\^eme, on a  
\[
\ttv(Q_m(u(x),y)-\sum_{j=0}^mx_j^{p^{m-j}}y_j)
\geq -n+2(r-1)\geq -n+r, 
\]
ce qui prouve \eqref{ver31}. 

Soit $x\in \rV^{m+1-m'}(\fil_{n+1}\rW_{m'}(K))\subset 
\fil_{n+1}\rW_{m+1}(K)$, o\`u $m'=\inf(\ord_p(n+1),m+1)$. On a
\[
p^{m-i}\ttv(Q_i(u(x),y)-\sum_{j=0}^ix_j^{p^{i-j}}y_j)
\geq -n-1+2p^{m-i}(r-1)\geq -n+r,
\]
ce qui d\'emontre \eqref{ver31}.
On dira que $i$ est ordinaire si $p^{m-i}\ttv(x_i)\geq -n$;  
$i$ est exceptionnel si $p^{m-i}\ttv(x_i)=-n-1$. 
Si $i$ est exceptionnel, alors $i\geq m+1-m'$; donc $\ord_p(\ttv(x_i))\geq 1$
et $\ttv(y_i)\geq r$ \eqref{vx3}.  
Soit $0\leq j\leq i$. Si $j$ est ordinaire, alors 
\[
p^{m-i}\ttv(x_j^{p^{i-j}}y_j)\geq -n+p^{m-i}(r-1). 
\]
Si $j$ est exceptionnel, alors
\[
p^{m-i}\ttv(x_j^{p^{i-j}}y_j)\geq -n-1+p^{m-i}r. 
\]
Les relations \eqref{ver11} et \eqref{ver21} s'en d\'eduisent aussit\^ot.

\begin{rema}
Il r\'esulte de la preuve de \ref{ver} que 
la condition $r\geq 2$ suffit 
sauf peut \^etre pour \eqref{ver31} si $x\not\in\fil_n\rW_{m+1}(K)$. 
\end{rema}
 
\begin{prop}\label{ig2}
Supposons $r=\ttv(D)\geq 3$. 
Pour tout entiers $m\geq 0$ et $n\geq r-1$, l'homomorphisme de groupes
\begin{equation}
v-u\colon \rW_{m+1}(K)\rightarrow  \rW_{m+1}(\mK)
\end{equation}
envoie $\fil'_n\rW_{m+1}(K)$ dans  $\fil_{n-r+1}\rW_{m+1}(\mK)$
et le diagramme 
\begin{equation}
\xymatrix{
{\Gr'_n\rW_{m+1}(K)}\ar[rr]^{\gr(v-u)}\ar[d]_{\gr'(\rF^md)}&&
{\Gr_{n-r+1}\rW_{m+1}(\mK)}\\
{\Gr'_n\Omega^1_K}\ar[rr]_{t_n}&&{\Gr_{n-r+1}\rW_1(\mK)}\ar[u]_{\gr(\rV^{m})}}
\end{equation}
est commutatif. 
\end{prop}

Cela r\'esulte de \ref{witt}(iii), \ref{ver} et \eqref{vx3}.

\begin{lem}\label{ig3}
Supposons $r=\ttv(D)\geq 2$. Soient $n\geq 1$ un entier premier \`a $p$,
$\chi\in \rH^1(K,\mu_n)$. Alors on a $v^*(\chi)=u^*(\chi)$ dans 
$\rH^1(\mK,\mu_n)$. 
\end{lem}

En vertu de \eqref{vx2}, pour tout $x\in K^\times$, 
$v(x)/u(x)\in 1+\fm_\mK$; c'est donc une puissance $n$-\`eme de $R_\mK$. 
La proposition r\'esulte alors de la th\'eorie de Kummer.

\begin{prop}\label{ig4}
Soient $n\geq 1$ un entier, $n'$ le plus grand diviseur de $n$ premier
\`a $p$. On suppose que $K$ contient une racine primitive $n'$-\`eme de
l'unit\'e. Soit $\chi\in \rH^1(K,\mZ/n\mZ)$ tel que 
$\sw'(\chi)\geq 2$ et $r=\ttv(D)=\sw'(\chi)+1$. 
Alors $n$ est un multiple de $p$ et 
$v^*(\chi)-u^*(\chi)$ est l'image de l'\'el\'ement 
$-t_{\sw'(\chi)}(\rsw'(\chi))$ par l'homomorphisme compos\'e 
\[
R_\mK/\fm_\mK \rightarrow \rH^1_\et(\Spec(R_\mK/\fm_\mK),\mZ/p\mZ)
\stackrel{\sim}{\rightarrow}
\rH^1_\et(\Spec(R_\mK),\mZ/p\mZ)\rightarrow \rH^1(\mK,\mZ/p\mZ)
\rightarrow\rH^1(\mK,\mZ/n\mZ),
\]
o\`u la premi\`ere fl\`eche 
est l'homomorphisme d'Artin-Schreier et la derni\`ere fl\`eche 
est induite par la multiplication par $n/p$ sur les coefficients.  
\end{prop}

Cela r\'esulte de \ref{ig2} et \ref{ig3}.

\subsection{}\label{ram-1}
On conserve dans la suite de cette section les hypoth\`eses de \eqref{ig4},  
de plus on fixe un homomorphisme injectif 
$\mZ/n\mZ\rightarrow \Lambda^\times$ tel que 
le compos\'e $\mZ/p\mZ \rightarrow 
\mZ/n\mZ\rightarrow \Lambda^\times$ avec la multiplication par $n/p$
soit le caract\`ere $\psi$ \eqref{fpsi}. 
On consid\`ere $-\rsw'(\chi)$ comme un $\oF$-point de $\obT_D^\ttt$, 
ou de fa\c{c}on \'equivalente comme une forme lin\'eaire 
$f\colon \obT_D\rightarrow \mA^1_{\os}$. 
On pose $\cL_\chi=f^*(\cL_\psi)$, o\`u $\cL_\psi$ est le 
faisceau d'Artin-Schreier associ\'e \`a $\psi$ (\cite{laumon} 1.1.3).  
On d\'esigne par $\cF$ le faisceau 
\'etale localement constant constructible en $\Lambda$-modules 
de rang $1$ sur $\eta$ associ\'e au caract\`ere $\chi$. 
En vertu de \ref{alg15}, quitte \`a remplacer $(X,\ttx)$ 
par un objet de $\fE$, 
il existe $\chi_U\in \rH^1(U,\mZ/n\mZ)$ tel que $\xi_U^*(\chi_U)=\chi$. 
Posons $\cG$ le faisceau \'etale localement constant constructible 
en $\Lambda$-modules de rang $1$ sur $U$ associ\'e au caract\`ere~$\chi_U$.

Comme $r=\ttv(D)\geq 3$, $X_0\times_k\eta$ est ferm\'e dans 
$(X\times_kS)_{(D)}$ \eqref{dilat18}. 
Soit $(X\times_kS)^\circ_{(D)}$ l'ouvert compl\'ementaire de 
$X_0\times_k\eta$ dans $(X\times_kS)_{(D)}$, de sorte qu'on ait 
un diagramme commutatif canonique \`a carr\'es cart\'esiens 
\[ 
\xymatrix{
{\bT_D}\ar[r]\ar[d]&{(X\times_kS)^{\circ}_{(D)}}\ar[d]&
{U\times_k\eta}\ar[l]\ar[d]&\\
D\ar[r]&{S}&{\eta}\ar[l]}
\]
On note $\Psi_D:\bD^+(U\times_k\eta)\rightarrow \bD^+(\obT_{D})$
le foncteur des cycles proches.

\begin{prop}\label{ram-2}
Le complexe $\Psi_D(\cH(\cF,\cG))$ \eqref{CD1} est 
isomorphe \`a $\cL_\chi[0]$.
\end{prop}
Le morphisme $(X\times_kS)^{\circ}_{(D)}\rightarrow S$
\'etant lisse, il suffit de montrer que $\cH(\cF,\cG)$ 
se prolonge en un faisceau 
localement constant constructible sur $(X\times_kS)^\circ_{(D)}$ 
dont la restriction \`a $\obT_D$ est isomorphe \`a $\cL_\chi$. 
Cela r\'esulte de \ref{ig4} par le th\'eor\`eme de puret\'e de 
Zariski-Nagata (SGA 2 X 3.4) et SGA 1 V 8.2.

\subsection{}\label{ram-3}
Soit $L$ la sous-extension de $\oK$ fix\'ee par le noyau du caract\`ere 
$\chi\colon G\rightarrow \mZ/n\mZ$. En vertu de \ref{alg15}, 
quitte \`a remplacer $(X,\ttx)$ par un objet de $\fE$, il existe 
un morphisme fini $f\colon Y\rightarrow X$ et un diagramme cart\'esien 
\begin{equation}\label{dram1}
\xymatrix{
{S_L}\ar[r]^\rho\ar[d]&Y\ar[d]^f\\
S\ar[r]_\xi&X}
\end{equation}
tels que $\rho$ induise un isomorphisme entre les groupes d'automorphismes
$\Aut_{X}(Y)\simeq \Aut_K(L)=\mZ/n\mZ$ et
$f$ soit \'etale au dessus de $U$  
de classe $\chi_U\in \rH^1(U,\mZ/n\mZ)$. 
Donc $Y$ a un unique 
point $\tty$ au dessus de $\ttx$, et son anneau local compl\'et\'e 
est isomorphe \`a $R_L$. 
On d\'esigne par $H_{(D)}$ le sch\'ema d\'efini par le diagramme cart\'esien 
\[
\xymatrix{
{H_{(D)}}\ar[r]\ar[d]_{f_{(D)}}&{Y\times_kS}\ar[d]^{f\times_kS}\\
{(X\times_kS)^\circ_{(D)}}\ar[r]&{X\times_kS}}
\]
par $H'_{(D)}$ la normalisation de $H_{(D)}\times_SS_L$
dans $Y_U\times_k\eta_L$ et par $f'_{(D)}\colon H'_{(D)}\rightarrow 
(X\times_kS)^\circ_{(D)}\times_SS_L$ le morphisme d\'eduit de $f_{(D)}$.

\begin{prop}\label{ram-4}
{\rm (i)}\ Le morphisme $f'_{(D)}$ est \'etale, et il induit un rev\^etement 
\'etale non-trivial de $\bT_D\times_SS_L$.

{\rm (ii)}\ On a $\chi(G^{r+})=0$ et  $\chi(G^{r})\not=0$, 
o\`u $r=\sw'(\chi)+1$.
\end{prop}
(i) Soit $\kappa'$ le point g\'en\'erique de $\bT_D\times_SS_L$. 
Le compl\'et\'e s\'epar\'e de l'anneau local de 
$(X\times_kS)_{(D)}\times_SS_L$ en $\kappa'$ s'identifie \`a 
$R_\mK\otimes_RR_L$. 
Il r\'esulte de \ref{ig4} que $f'_{(D)}$ est \'etale en tout point au dessus  
de $\kappa'$ et le rev\^etement \'etale induit par 
$f'_{(D)}$ au dessus de $\kappa'$ est non-trivial (\ref{msw3} et \ref{msw35}).
On en d\'eduit par le th\'eor\`eme de puret\'e de Zariski-Nagata 
que $f'_{(D)}$ est un rev\^etement \'etale. 

(ii) On sous-entend par compl\'et\'e formel d'un $S$-sch\'ema 
son compl\'et\'e formel le long de sa fibre sp\'eciale. 
Consid\'erons le diagramme cart\'esien 
\begin{equation}\label{dram2}
\xymatrix{
{S_L}\ar[r]^-(0.5){\theta}\ar[d]&{Y\times_kS}\ar[d]^{f\times_kS}\\
S\ar[r]_-(0.5){\gamma}&{X\times_kS}}
\end{equation}
induit par \eqref{dram1}, o\`u $\gamma$ est le graphe de $\xi$. 
On note $\fX$ (resp. $\fY$) le compl\'et\'e 
formel de $X\times_kS$ (resp. $Y\times_kS$) 
et $\fX^\rig$ (resp. $\fY^\rig$) sa fibre rigide au sens de Raynaud. 
Pour tout nombre rationnel $j\geq 0$, 
on d\'esigne par $\fX^{\rig}_{(j)}$ (resp. $\fY^{\rig}_{(j)}$) le
voisinage tubulaire de $\gamma$ (resp. $\theta$) de rayon $j$ 
(\cite{as} part II 1.2); 
c'est un sous-domaine affino\"{\i}de de $\fX^\rig$ (resp. $\fY^\rig$).
On a $\fY^\rig_{(j)}=\fX^\rig_{(j)}\times_{\fX^\rig}\fY^\rig$; 
on note $f^\rig_{(j)}\colon \fY^\rig_{(j)}\rightarrow \fX^\rig_{(j)}$
la projection canonique. Le compl\'et\'e formel de $(X\times_kS)^\circ_{(D)}$
est un mod\`ele formel de $\fX^{\rig}_{(r)}$. 
Les morphismes $f_{(D)}$ et $f'_{(D)}$ \'etant finis, 
leurs compl\'et\'es formels sont alors des mod\`eles formels 
des morphismes $f_{(r)}^\rig$ et $f_{(r)}^\rig\otimes_KL$ respectivement. 
La filtration $(G^j)_{j\in \mQ_{\geq 0}}$ de $G$ est d\'efinie dans 
\cite{as} de sorte que 
\[
(\mZ/n\mZ)/\chi(G^j)\simeq \pi_0(\fY^\rig_{(j),\oK}),
\]
o\`u le terme de droite d\'esigne l'ensemble des composantes connexes 
g\'eom\'etriques de $\fY^\rig_{(j)}$. 
Par cons\'equent, la proposition (ii) r\'esulte de (i) 
(\cite{as} part I, 4.3).

\subsection{}\label{ptp1} 
Venons maintenant \`a la preuve du th\'eor\`eme \ref{tp1}. 
Soit $\ok$ une cl\^oture alg\'ebrique de $k$. Il suffit de d\'emontrer 
\ref{tp1} apr\`es avoir remplac\'e $S$ par une composante connexe 
de $S_{\ok}$. Donc on peut supposer que $K$ contient une racine 
primitive $n'$-\`eme de l'unit\'e, o\`u $n'$ est le plus grand diviseur 
de $n$ premier \`a $p$. 
Il r\'esulte alors de \ref{ram-4}(ii) que $r=\sw'(\chi)+1$ est l'unique 
pente critique de $\cF$. Compte tenu de la d\'emonstration de \ref{CD2},   
on a un isomorphisme canonique $\nu_D(\cH(\cF_\oeta,\cF_!))\simeq 
\Psi_D(\cH(\cF,\cG))$; donc en vertu de \ref{ram-2}, on a 
$\nu_D(\cH(\cF_\oeta,\cF_!))\simeq \cL_\chi[0]$.  
Par ailleurs, il r\'esulte par exemple de \cite{laumon}
1.2.3.2 et 1.2.2.1, que la transform\'ee de Fourier du faisceau 
$\cL_\chi=f^*(\cL_\psi)$ est un faisceau de support la section 
de $\obT^{\ttt}_D$ associ\'ee au morphisme 
$f\colon \obT_D\rightarrow \mA^1_k$ (c'est \`a dire la section $-\rsw'(\chi)$),
d'o\`u l'assertion \ref{tp1}(ii).

\section{D\'emonstration du th\'eor\`eme \ref{tp2}}\label{Dtp2}

\subsection{}\label{jg1}
Dans cette section, on d\'esigne par $(X,\ttx)$ le $k$-sch\'ema point\'e 
et $\xi\colon S\rightarrow X$ le $k$-morphisme 
du \eqref{alg1}, par $D$ un diviseur effectif de $S$. 
On note $X_0$ l'adh\'erence sch\'ematique de $\ttx$ dans $X$
(qui est un diviseur de Cartier) et  
$U$ l'ouvert compl\'ementaire de $X_0$ dans $X$. On pose 
$X\times_k^{\log}S$ le produit fibr\'e logarithmique
relativement aux diviseurs de Cartier $X_0$ sur $X$ et $s$ sur $S$ 
\eqref{plog3}, $r=\ttv(D)$ \eqref{not2},   
$\Theta_D=\bV(\Omega^1_R(\log)\otimes_R\co_D(D))$
et $\oTheta_D=\Theta_D\times_D\os$.
Comme $\xi^{-1}(X_0)=s$, on peut appliquer 
la construction \eqref{plog4} au morphisme
$\xi$; on d\'esigne par $(X\times_k^{\log}S)_{(D)}$ la dilatation de
$X\times_k^{\log}S$ le long du graphe de $\xi$ d'\'epaisseur $D$. 
On a alors un diagramme commutatif canonique \`a carr\'es cart\'esiens 
\begin{equation} 
\xymatrix{
{\Theta_D}\ar[r]\ar[d]&{(X\times_k^{\log}S)_{(D)}}\ar[d]&
{U\times_k\eta}\ar[l]\ar[d]&\\
D\ar[r]&{S}&{\eta}\ar[l]}
\end{equation}
On note $\kappa$ le point g\'en\'erique de $\Theta_D$, 
$R_\mK$ le s\'epar\'e compl\'et\'e de l'anneau local 
de $(X\times_k^{\log}S)_{(D)}$ en $\kappa$, $\mK$ son corps des fractions, 
$\fm_\mK$ son id\'eal maximal. 
La projection canonique $(X\times_k^{\log}S)_{(D)}\rightarrow S$ \'etant 
lisse, elle induit un homomorphisme $u\colon R\rightarrow R_\mK$
v\'erifiant $\fm_\mK=u(\fm)R_\mK$. On note encore
$\ttv\colon \mK^\times \rightarrow \mZ$ la valuation qui prolonge 
$\ttv$ sur $K$ \eqref{not1}.  
Le point $\kappa$ est uniquement d\'etermin\'e par un $F$-homomorphisme 
\begin{equation}\label{lder-t1}
t\colon \Omega^1_R(\log)\rightarrow \fm_\mK^{r}/\fm_\mK^{r+1}.
\end{equation}
Appliquant le produit tensoriel $\otimes_F(\fm^{-n}/\fm^{-n+1})$ 
(pour $n\in \mZ$), on obtient un $F$-homomorphisme
\begin{equation}
t_n\colon \Gr_n\Omega^1_K\rightarrow \fm_\mK^{-n+r}/\fm_\mK^{-n+r+1}
=\Gr_{n-r}\rW_1(\mK).
\end{equation}
La premi\`ere projection $(X\times_k^{\log}S)_{(D)}\rightarrow X$ 
envoie $\kappa$ sur $\ttx$. On en d\'eduit un second homomorphisme 
$v\colon R\rightarrow R_\mK$.

\begin{lem}\label{lvx}
Pour tout $x\in K^\times$, on a $\ttv(v(x)/u(x)-1)\geq r$ et 
\[
\frac{v(x)}{u(x)}-1\equiv t(d\log(x))\mod \fm_\mK^{r+1}.
\]
\end{lem}
Compte tenu de la relation \eqref{vx4}, on peut se borner aux cas o\`u 
$x$ est une unit\'e de $R$ ou $x=\pi$, car $r\geq 1$. 
Le premier cas est \'evident (cf. \ref{vx}(i)). Le second
se d\'eduit d'une relation 
analogue et \'evidente dans l'anneau 
$R\otimes_kR[\ttu,\ttu^{-1}]/(\pi\otimes 1-1\otimes \pi \cdot \ttu)$.

\begin{lem} \label{lver}
Soient $m,n\geq 0$ deux entiers, $x=(x_0,\dots,x_m)\in \fil_{n}\rW_{m+1}(K)$, 
$y=(y_0,\dots,y_m)\in \mK^{m+1}$, o\`u 
$y_i=\frac{v(x_i)}{u(x_i)}-1$ si $ x_i\not=0$ et 
$y_i=0$ si $x_i=0$. Alors on a
\begin{eqnarray}
p^{m-i}\ttv(Q_i(u(x),y))&\geq& -n+r, \ \ \ \ \ \ \ \forall 
\ 0\leq i\leq m,\label{lver11}\\
\ttv(Q_m(u(x),y)- \sum_{j=0}^mx_j^{p^{m-j}}y_j)&\geq& -n+2r\geq -n+r+1.
\label{lver31}
\end{eqnarray}
\end{lem}

Cela r\'esulte de \ref{lvx} en calquant la preuve de \ref{ver}.

\begin{prop}\label{jg2}
Pour tout entiers $m\geq 0$ et $n\geq r$, l'homomorphisme de groupes
\begin{equation}
v-u\colon \rW_{m+1}(K)\rightarrow  \rW_{m+1}(\mK)
\end{equation}
envoie $\fil_n\rW_{m+1}(K)$ dans  $\fil_{n-r}\rW_{m+1}(\mK)$
et le diagramme 
\begin{equation}
\xymatrix{
{\Gr_n\rW_{m+1}(K)}\ar[rr]^{\gr(v-u)}\ar[d]_{\gr(\rF^md)}&&
{\Gr_{n-r}\rW_{m+1}(\mK)}\\
{\Gr_n\Omega^1_K}\ar[rr]_{t_n}&&{\Gr_{n-r}\rW_1(\mK)}\ar[u]_{\gr(\rV^{m})}}
\end{equation}
est commutatif. 
\end{prop}

Cela r\'esulte de \ref{witt}(iii), \ref{lvx} et \ref{lver}.

\begin{lem}\label{jg3}
Soient $n\geq 1$ un entier premier \`a $p$,
$\chi\in \rH^1(K,\mu_n)$. Alors on a $v^*(\chi)=u^*(\chi)$ dans 
$\rH^1(\mK,\mu_n)$. 
\end{lem}

En vertu de \ref{lvx}, pour tout $x\in K^\times$, 
$v(x)/u(x)\in 1+\fm_\mK$; c'est donc une puissance $n$-\`eme de $R_\mK$. 
La proposition r\'esulte alors de la th\'eorie de Kummer.

\begin{prop}\label{jg4}
Soient $n\geq 1$ un entier, $n'$ le plus grand diviseur de $n$ premier
\`a $p$. On suppose que $K$ contient une racine primitive $n'$-\`eme de
l'unit\'e. Soit $\chi\in \rH^1(K,\mZ/n\mZ)$ tel que 
$\sw(\chi)\geq 1$ et $r=\ttv(D)=\sw(\chi)$. 
Alors $n$ est un multiple de $p$ et 
$v^*(\chi)-u^*(\chi)$ est l'image de l'\'el\'ement 
$-t_{\sw(\chi)}(\rsw(\chi))$ par l'homomorphisme compos\'e 
\[
R_\mK/\fm_\mK \rightarrow \rH^1_\et(\Spec(R_\mK/\fm_\mK),\mZ/p\mZ)
\stackrel{\sim}{\rightarrow}
\rH^1_\et(\Spec(R_\mK),\mZ/p\mZ)\rightarrow \rH^1(\mK,\mZ/p\mZ)
\rightarrow\rH^1(\mK,\mZ/n\mZ),
\]
o\`u la premi\`ere fl\`eche 
est l'homomorphisme d'Artin-Schreier et la derni\`ere fl\`eche 
est induite par la multiplication par $n/p$ sur les coefficients.  
\end{prop}
Cela r\'esulte de \ref{jg2} et \ref{jg3}. 

\subsection{}\label{lram-1}
On conserve dans la suite de cette section les hypoth\`eses de \eqref{jg4}, 
de plus on fixe un homomorphisme $\mZ/n\mZ\rightarrow \Lambda^\times$ tel que 
le compos\'e $\mZ/p\mZ \rightarrow 
\mZ/n\mZ\rightarrow \Lambda^\times$ avec la multiplication par $n/p$
soit le caract\`ere $\psi$ \eqref{fpsi}. 
On consid\`ere $-\rsw(\chi)$ comme un $\oF$-point de $\oTheta_D^\ttt$, 
ou de fa\c{c}on \'equivalente comme une forme lin\'eaire 
$f\colon \oTheta_D\rightarrow \mA^1_{\os}$. 
On pose $\cL_\chi=f^*(\cL_\psi)$, o\`u $\cL_\psi$ est le 
faisceau d'Artin-Schreier associ\'e \`a $\psi$ (\cite{laumon} 1.1.3).  
On d\'esigne par $\cF$ le faisceau 
\'etale localement constant constructible en $\Lambda$-modules 
de rang $1$ sur $\eta$ associ\'e au caract\`ere $\chi$. 
En vertu de \ref{alg15}, quitte \`a remplacer $(X,\ttx)$ 
par un objet de $\fE$, 
il existe $\chi_U\in \rH^1(U,\mZ/n\mZ)$ tel que $\xi_U^*(\chi_U)=\chi$. 
Posons $\cG$ le faisceau \'etale localement constant constructible 
en $\Lambda$-modules de rang $1$ sur $U$ associ\'e au caract\`ere $\chi_U$.
On note $\Psi_D:\bD^+(U\times_k\eta)\rightarrow \bD^+(\oTheta_{D})$
le foncteur des cycles proches pour le morphisme $(X\times_k^{\log}S)_{(D)}
\rightarrow S$.

\begin{prop}\label{lram-2}
Le complexe $\Psi_D(\cH(\cF,\cG))$ \eqref{CD1} est 
isomorphe \`a $\cL_\chi[0]$.
\end{prop}
Le morphisme $(X\times_k^{\log}S)_{(D)}
\rightarrow S$ \'etant lisse, 
il suffit de montrer que $\cH(\cF,\cG)$ se prolonge en un faisceau 
localement constant constructible sur $(X\times_k^{\log}S)_{(D)}$ 
dont la restriction 
\`a $\oTheta_D$ est isomorphe \`a $\cL_\chi$. 
Cela r\'esulte de \ref{jg4} par le th\'eor\`eme de puret\'e de 
Zariski-Nagata (SGA 2 X 3.4) et SGA 1 V 8.2.

\subsection{}\label{lram-3}
Soit $L$ la sous-extension de $\oK$ fix\'ee par le noyau du caract\`ere 
$\chi\colon G\rightarrow \mZ/n\mZ$. En vertu de \ref{alg2},  
quitte \`a remplacer $(X,\ttx)$ par un objet de $\fE$, il existe 
un morphisme fini $f\colon Y\rightarrow X$ et un diagramme cart\'esien 
\begin{equation}\label{ldram1}
\xymatrix{
{S_L}\ar[r]^\rho\ar[d]&Y\ar[d]^f\\
S\ar[r]_\xi&X}
\end{equation}
tels que $\rho$ induise un isomorphisme entre les groupes d'automorphismes
$\Aut_{X}(Y)\simeq \Aut_K(L)=\mZ/n\mZ$ et
$f$ soit \'etale au dessus de $U$  
de classe $\chi_U\in \rH^1(U,\mZ/n\mZ)$.  
On d\'esigne par $H$ et $H_{(D)}$ les sch\'emas d\'efinis 
par le diagramme cart\'esien 
\[
\xymatrix{
{H_{(D)}}\ar[r]\ar[d]_{h_{(D)}}&{H}\ar[d]^{h}\ar[r]&Y\ar[d]^f\\
{(X\times_k^{\log}S)_{(D)}}\ar[r]&{(X\times_k^{\log}S)}\ar[r]&X}
\]
par $H'_{(D)}$ la normalisation de $H_{(D)}\times_SS_L$
dans $Y_U\times_k\eta_L$ et par $h'_{(D)}\colon H'_{(D)}\rightarrow 
(X\times_k^{\log}S)_{(D)}\times_SS_L$ le morphisme d\'eduit de $h_{(D)}$.

\begin{prop}\label{lram-4}
{\rm (i)}\ Le morphisme $h'_{(D)}$ est \'etale, et il induit un rev\^etement 
\'etale non-trivial de $\Theta_D\times_SS_L$.

{\rm (ii)}\ On a $\chi(G^{r+}_{\log})=0$ et  $\chi(G^{r}_{\log})\not=0$, 
o\`u $r=\sw(\chi)$
\end{prop}
(i) Soit $\kappa'$ le point g\'en\'erique de $\Theta_D\times_SS_L$. 
Le compl\'et\'e s\'epar\'e de l'anneau local de 
$(X\times_k^{\log}S)_{(D)}\times_SS_L$ en $\kappa'$ s'identifie \`a 
$R_\mK\otimes_RR_L$. Il r\'esulte de \ref{jg4} que $h'_{(D)}$ est \'etale 
en tout point au dessus  de $\kappa'$ et le rev\^etement \'etale induit par 
$h'_{(D)}$ au dessus de $\kappa'$ est non-trivial (\ref{sw3} et \ref{sw35}).
On en d\'eduit par le th\'eor\`eme de puret\'e de Zariski-Nagata que 
$h'_{(D)}$ est un rev\^etement \'etale.

(ii) On sous-entend par compl\'et\'e formel d'un $S$-sch\'ema 
son compl\'et\'e formel le long de sa fibre sp\'eciale. 
Consid\'erons le diagramme cart\'esien 
\begin{equation}\label{ldram2}
\xymatrix{
{S_L}\ar[r]^-(0.5){\theta}\ar[d]&{H}\ar[d]^{h}\\
S\ar[r]_-(0.5){\gamma^{\log}}&{(X\times_k^{\log}S)}}
\end{equation} 
induit par \eqref{ldram1}, o\`u 
$\gamma^{\log}$ est le graphe logarithmique de $\xi$ \eqref{plog3}. 
On note $\fX$ (resp. $\fY$) le compl\'et\'e formel de $(X\times_k^{\log}S)$ 
(resp. $H$) et $\fX^\rig$ (resp. $\fY^\rig$) sa fibre rigide. 
Pour tout nombre rationnel $j\geq 0$, 
on d\'esigne par $\fX^{\rig}_{(j)}$ (resp. $\fY^{\rig}_{(j)}$) le
voisinage tubulaire de $\gamma^{\log}$ (resp. $\theta$) de rayon $j$ 
(\cite{as} II 1.2). On a 
$\fY^\rig_{(j)}=\fX^\rig_{(j)}\times_{\fX^\rig}\fY^\rig$; 
on note $h^\rig_{(j)}\colon \fY^\rig_{(j)}\rightarrow \fX^\rig_{(j)}$
la projection canonique. Le compl\'et\'e formel de 
$(X\times_k^{\log}S)_{(D)}$ est un mod\`ele formel de $\fX^{\rig}_{(r)}$.
Les morphismes $h_{(D)}$ et $h'_{(D)}$ \'etant finis, 
leurs compl\'et\'es formels sont alors des mod\`eles formels 
des morphismes $h_{(r)}^\rig$ et $h_{(r)}^\rig\otimes_KL$ respectivement. 
La filtration $(G^j_{\log})_{j\in \mQ_{\geq 0}}$ de 
$G$ est d\'efinie dans \cite{as} de sorte que 
\[
(\mZ/n\mZ)/\chi(G^j_{\log})\simeq \pi_0(\fY^\rig_{(j),\oK}),
\]
o\`u le terme de droite d\'esigne l'ensemble des composantes connexes 
g\'eom\'etriques de $\fY^\rig_{(j)}$. 
Par cons\'equent, la proposition (ii) r\'esulte de (i) (\cite{as} part I, 4.3).

\subsection{}\label{ptp2} 
On peut maintenant d\'emontrer le th\'eor\`eme \ref{tp2} 
en calquant la preuve de \ref{tp1} (cf. \ref{ptp1}).

\section{Ramification des corps locaux \`a corps r\'esiduel parfait}
\label{Dparf}

\subsection{}\label{hypf}
On suppose que $F=k$ et on se propose de calculer pour une 
$\Lambda$-repr\'esentation finie de $G$ 
les invariants de ramification d\'efinis dans cet article 
en termes d'invariants plus classiques.
Pour tout nombre rationnel $r>0$, on sait (\cite{as} part I 3.7 et 3.15)
que $G^{r+1}=G_{\log}^r$ est le $r$-\`eme sous-groupe 
de ramification sup\'erieure ``classique''de $G$ (\cite{serre1} IV Section 3).

\subsection{}\label{id-can}
Pour tout nombre rationnel $r$, on pose $\rN_r=\fm^{(r)}
\otimes_{\oR}\oF$ (cf. \ref{mj}); c'est un $\oF$-espace vectoriel de 
dimension $1$ dont le dual est canoniquement isomorphe \`a $\rN_{-r}$; 
donc $\bV(\rN_r)(\oF)=\rN_{-r}$. On observe que si $r>0$, 
$\obT_r$ s'identifie canoniquement \`a 
$\bV(\rN_{1-r})$, et $\oTheta_r$ \`a $\bV(\rN_{-r})$ gr\^ace \`a \eqref{res1}.

\subsection{} 
Pour tout nombre rationnel $r$,
on d\'esigne par $\pi_1^{\iso}(\bV(\rN_r))$ le quotient du groupe 
fondamental de $\bV(\rN_r)$ qui classifie les isog\'enies \'etales 
de groupes alg\'ebriques; c'est un groupe profini, ab\'elien et
annul\'e par $p$.
On rappelle que l'isog\'enie de Lang $\mA^1_{\os}\rightarrow 
\mA^1_{\os}$, d\'efinie par $x\mapsto x^p-x$, est une base du 
$\oF$-vectoriel $\Hom_\mZ(\pi_1^{\iso}(\mA^1_{\os}),\mF_p)$. 
On en d\'eduit un isomorphisme canonique 
\begin{equation}\label{iso-dual}
\Hom_\mZ(\pi_1^{\iso}(\bV(\rN_r)),\mF_p)\simeq 
\bV(\rN_{-r})(\oF)=\rN_r.
\end{equation}
L'action de $G$ sur $\oK$ induit une action sur $\rN_r$
et donc sur $\pi_1^\iso(\bV(\rN_r))$.

\begin{teo}[\cite{as} part II 6.3(2)]\label{grad63}
Pour tout nombre rationnel $r>0$, 
on a un isomorphisme canonique $G$-\'equivariant 
\begin{equation}\label{grad64}
\pi_1^{\iso}(\bV(\rN_{-r}))\simeq G^{r}_{\log}/G^{r+}_{\log},
\end{equation}
o\`u l'action de $G$ sur $G^{r}_{\log}/G^{r+}_{\log}$ est induite 
par la conjugaison. 
\end{teo}

\subsection{}\label{ram-cons}
Soient $L$ une extension galoisienne finie de $K$ de groupe $\oG$,
$(\oG^j_{\log})_{j\in \mQ_{\geq 0}}$ la filtration de ramification 
sup\'erieure logarithmique de $\oG$, 
$r$ le {\em conducteur logarithmique} de $L$ sur $K$,
c'est \`a dire l'unique nombre rationnel $r\geq 0$ tel que 
$\oG^r_{\log}\not=0$ et $\oG^{r+}_{\log}=0$.
On suppose que $L/K$ est sauvagement ramifi\'ee, de sorte que $r>0$. 
On sait qu'il existe un diviseur $D$ de $S_L$ tel que $\ttv(D)=r$
\eqref{not2}. 

On d\'esigne par $(X,\ttx)$ le $k$-sch\'ema point\'e 
et $\xi\colon S\rightarrow X$ le $k$-morphisme 
du \eqref{alg1}, par $X_0$ l'adh\'erence sch\'ematique de $\ttx$ dans $X$
(qui est un diviseur de Cartier) et 
par $U$ l'ouvert compl\'emen\-taire de $X_0$ dans $X$. 
On pose $q\colon S_L\rightarrow S$ la fl\`eche canonique,
$X\times_k^{\log}S_L$ le produit fibr\'e logarithmique
relativement aux diviseurs de Cartier $X_0$ sur $X$ et $q^{-1}(s)$ sur $S_L$ 
\eqref{plog3}. Comme $\xi^{-1}(X_0)=s$, on peut appliquer 
la construction \eqref{plog4} au morphisme
$\xi\circ q$; on d\'esigne par $(X\times_k^{\log}S_L)_{(D)}$ la dilatation de
$X\times_k^{\log}S_L$ le long du graphe de $\xi\circ q$ d'\'epaisseur $D$. 
On a alors un diagramme commutatif canonique \`a carr\'es cart\'esiens 
\begin{equation} 
\xymatrix{
{\Theta_D}\ar[r]\ar[d]&{(X\times_k^{\log}S_L)_{(D)}}\ar[d]&
{U\times_k\eta_L}\ar[l]\ar[d]&\\
D\ar[r]&{S_L}&{\eta_L}\ar[l]}
\end{equation}
o\`u $\Theta_D=\bV(\Omega^1_R(\log)\otimes_R\co_D(D))$. 

En vertu de \ref{alg15}, quitte \`a remplacer $(X,\ttx)$ par un objet 
de $\fE$, il existe 
un morphisme fini $f\colon Y\rightarrow X$ et un diagramme cart\'esien 
\begin{equation}\label{hq1}
\xymatrix{
{S_L}\ar[r]^\rho\ar[d]&Y\ar[d]^f\\
S\ar[r]_\xi&X}
\end{equation}
tels que $\rho$ induise un isomorphisme entre les groupes d'automorphismes
$\Aut_{X}(Y)\simeq \Aut_K(L)=\oG$ et $f$ soit \'etale au dessus de $U$. 
On d\'esigne par $H_{(D)}$ la normalisation de $(X\times_k^{\log}S_L)_{(D)}$
dans $Y_U\times_k\eta_L$ et par $f_{(D)}\colon H_{(D)}\rightarrow 
(X\times_k^{\log}S_L)_{(D)}$ le morphisme canonique. 

\begin{prop}\label{hq3}
Le morphisme $f_{(D)}$ est \'etale; 
il induit sur chaque composante connexe 
au dessus de $\oTheta_D=\bV(\rN_{-r})$ une isog\'enie \'etale de groupe 
$\oG^r_{\log}$; l'homomorphisme 
$\pi_1^\iso(\bV(\rN_{-r}))\rightarrow \oG^r_{\log}$ ainsi d\'efini est 
compatible \`a \eqref{grad64}. 
\end{prop}

Cela r\'esulte de \cite{as} part II 5.9(3), 6.3(2) et 5.10. 

\subsection{}\label{pparf2}
Venons maintenant \`a la preuve du th\'eor\`eme \ref{parf2}.
On se limite \`a d\'emontrer l'\'enonc\'e logarithmique \eqref{ciso2},
l'\'enonc\'e \eqref{ciso1} \'etant \'equivalent (la preuve de cet \'equivalence
est similaire \`a \cite{as} part I 9.11(ii) ou part II 6.4.1).   
Soient $\oG$ l'image de l'homomorphisme $\chi \colon 
G\rightarrow \Aut_\Lambda(\cF_\oeta)$,
$L$ la sous-extension de $\oK$ fix\'ee par le noyau de $\chi$. 
On reprend pour cet extension la construction \eqref{ram-cons}. 
Posons $\cG$ le faisceau \'etale localement constant constructible 
en $\Lambda$-modules sur $U$ associ\'e \`a la repr\'esentation 
$\chi$ (via le rev\^etement \'etale $f_U$). On note 
$\Psi_D:\bD^+(U\times_k\eta_L)\rightarrow \bD^+(\oTheta_{D})$
le foncteur des cycles proches pour le morphisme $(X\times_k^{\log}S_L)_{(D)}
\rightarrow S_L$. On a alors un isomorphisme canonique 
$\nu^{\log}_D(\cF\boxtimes\Lambda)\simeq \Psi_D(\cG\boxtimes\Lambda)$
(cf. la preuve de \ref{CD2}). 
D'apr\`es \ref{hq3}, $\cG\boxtimes\Lambda$ se prolonge en un faisceau 
localement constant 
constructible sur $(X\times_k^{\log}S_L)_{(D)}$. 
D'autre part, $(X\times_k^{\log}S_L)_{(D)}\rightarrow S_L$ 
est lisse. La conjecture \ref{ciso2} r\'esulte alors de \ref{hq3}.  

\subsection{}
On rappelle que $\Spec(\Lambda)$ est connexe \eqref{fpsi}.   
Soient $\cF$ un faisceau \'etale constructible en $\Lambda$-modules
sur $\eta$, $r>0$ une pente logarithmique critique de $\cF$ \eqref{pc}. 
On rappelle que $G^r_{\log}/G^{r+}_{\log}$ 
est un $\mF_p$-vectoriel et qu'on a fix\'e un caract\`ere non-trivial 
$\psi\colon \mF_p\rightarrow \Lambda^\times$. 
Donc l'action de $G^r_{\log}/G^{r+}_{\log}$ sur 
$(\cF_\oeta)^{G^{r+}_{\log}}$ d\'etermine un ensemble fini de caract\`eres 
$G^r_{\log}/G^{r+}_{\log}\rightarrow \mF_p$, 
ce qui induit par \eqref{iso-dual} 
et \eqref{grad64} un sous-ensemble fini $\Sigma_r(\cF)$ de $\rN_{-r}$. 

\begin{prop}\label{parf1}
Soient $\cF$ un faisceau \'etale constructible 
en $\Lambda$-modules plats sur $\eta$, 
$r$ sa plus grande pente logarithmique critique. 
Supposons que $\cF$ soit sauvagement ramifi\'e, de sorte que $r>0$.
Alors on a $\Sigma_r(\cF)=\ttC_{r+1}(\cF)=\trC_r(\cF)$ \eqref{id-can}.
\end{prop}

Il suffit de calquer \eqref{pparf2}.

\end{document}